\documentclass[reqno,12pt]{amsart}
\usepackage{latexsym,amssymb,amsmath,amsthm,bbm}
\usepackage{mathtools, color}
\usepackage{mathabx}

\newcommand{\R}{{\ensuremath{\mathbb{R}}}}
\newcommand{\N}{{\ensuremath{\mathbb{N}}}}
\newcommand{\Z}{{\ensuremath{\mathbb{Z}}}}

\newcommand{\B}{{\ensuremath{\mathcal{B}}}}
\renewcommand{\P}{\ensuremath{\mathbb{P}}}
\renewcommand{\dj}{d\kern-0.4em\char"16\kern-0.1em}
\newcommand{\E}{\ensuremath{\mathbb{E}}}

\newcommand{\wh}[1]{{\widehat{#1}}}

\newcommand{\cp}{\ensuremath{\mathrm{Cap}}}

\def\1{{\bf 1}}

\newtheorem{Thm}{Theorem}[section]
\newtheorem{Cor}[Thm]{Corollary}
\newtheorem{Lem}[Thm]{Lemma}
\newtheorem{Prop}[Thm]{Proposition}

\theoremstyle{definition}
\newtheorem{Rem}[Thm]{Remark}

\theoremstyle{definition}

\theoremstyle{definition}
\newtheorem{Def}[Thm]{Definition}

\theoremstyle{definition}

\def\sP {{\mathcal P}}

\def\wt{\widetilde}

\begin{document}
\numberwithin{equation}{section}

\author{Ante Mimica}
\address{Department of Mathematics, University of Zagreb,  Bijeni\v{c}ka c.~30, 10000 Zagreb, Croatia}
\curraddr{}
\thanks{}
\email{amimica@math.hr}
\author{Zoran Vondra\v{c}ek}
\address{Department of Mathematics, University of Zagreb, Bijeni\v{c}ka c.~30, 10000 Zagreb, Croatia}
\curraddr{}
\thanks{Supported in part by the MZOS grant 037-0372790-2801.}
\email{vondra@math.hr}

\begin{abstract}
A collection $\{\overline{B}(x_n,r_n)\}_{n\ge 1}$ of pairwise disjoint balls in the Euclidean space $\R^d$ is said to be avoidable with respect to a transient process $X$ if the process with positive probability escapes to infinity without hitting any ball. In this paper we study sufficient and necessary conditions for avoidability with respect to unimodal isotropic L\'evy processes satisfying a certain scaling hypothesis. These conditions are expressed in terms of the characteristic exponent of the process, or alternatively, in terms of the corresponding Green function. We also discuss the same problem for a random collection of balls. The results are generalization of several recent results for the case of Brownian motion.
\end{abstract}

\subjclass[2010]
{Primary 60J45, 60G51, Secondary 31B15 }

\keywords{Isotropic L\'evy process, Green function, minimal thinness at infinity}

\title{Unavoidable collections of balls for isotropic L\' evy processes}

\maketitle

\allowdisplaybreaks[3]

\section{Introduction}
Let $\{\overline{B}(x_n,r_n)\}_{n\ge 1}$ be a collection of pairwise disjoint balls in the Euclidean space $\R^d$, $d\ge 3$. For convenience we assume that none of the balls contains the origin. This collection of balls is said to be \emph{avoidable} (with respect to Brownian motion) if Brownian motion started at the origin with positive probability escapes to infinity without hitting any of the balls. More precisely, if $X=(X_t,\P_x)$ denotes a standard Brownian motion in $\R^d$,  $A=\bigcup_{n\ge 1}\overline{B}(x_n,r_n)$ is the union of all balls, and $T_A=\inf\{t>0:\ X_t\in A\}$  the hitting time to $A$, then  $A$ is avoidable if $\P_0(T_A < \infty)<1$. If the last probability is equal to 1, then we say that the collection of balls is unavoidable. The restriction on dimension $d\ge 3$ comes from the fact that for $d=1$ and $d=2$ any single ball is unavoidable.

The problem of determining when the collection of closed balls is avoidable has been studied in \cite{COC} and \cite{GG}. In particular, the following result was proved  in \cite{GG}.

\noindent
{\bf Theorem A.} \emph{If $A$ is unavoidable, then
\begin{equation}\label{e:thm-A-1}
\sum_{n\ge 1}\left(\frac{r_n}{|x_n|}\right)^{d-2}=\infty\, .
\end{equation}
Conversely, if \eqref{e:thm-A-1} and the separation condition
\begin{equation}\label{e:thm-A-2}
\inf_{m\neq n}\frac{|x_m-x_n|}{r_n^{1-2/d}|x_n|^{2/d}}>0
\end{equation}
hold, then $A$ is unavoidable.}

The same result with a stronger separation condition was proved in \cite[Theorem 1]{COC}.
It is shown in \cite[Theorem 3]{COC} that a separation condition is not redundant -- there exists an avoidable collection of balls $\{\overline{B}(x_n,r_n)\}_{n\ge 1}$ (satisfying  $\inf_{m\neq n}|x_m-x_n| >0$) such that \eqref{e:thm-A-1} holds.

The family of balls $\{\overline{B}(x_n,r_n)\}_{n\ge 1}$ is said to be \emph{regularly located} if the following three conditions are satisfied:
\begin{itemize}
    \item[(i)] There exists $\epsilon >0$ such that $|x_m-x_n|>2\epsilon$ for $m\neq n$ ($\epsilon$ separation);
    \item[(ii)] There exists $R>0$ such that any ball $B(x,R)$ contains at least one center $x_n$ (uniform density);
    \item[(iii)] There is a decreasing function $\phi:(0,\infty)\to (0,\infty)$ such that $r_n=\phi(|x_n|)$ for each $n\ge 1$ (the radius of the ball is the function of the distance of its center from the origin).
\end{itemize}
The following sufficient and necessary condition for avoidability of a regularly located collection of balls is proved in \cite[Theorem 2]{COC}.

\noindent
{\bf Theorem B.} \emph{Suppose that the collection of balls $\{\overline{B}(x_n,r_n)\}_{n\ge 1}$ is regularly located. Then this collection is avoidable if and only if
\begin{equation}\label{e:thm-B}
\int_1^{\infty} r\phi(r)^{d-2}\, dr <\infty\, .
\end{equation}
}

The main goal of this paper is to extend these two theorems to a class of isotropic L\'evy processes which, as a very special case, includes isotropic $\alpha$-stable processes. The extension is possible due to some very recent results in the potential theory of such L\'evy processes, cf.~\cite{G, KSV8}. We first note that since a L\'evy process can be transient in dimensions $d=1$ and $d=2$, the problem of avoidability can be studied in all dimensions. Still, throughout most of the paper we will focus on the case $d\ge 3$. The reason  is that in this case we need to impose minimal assumptions on the L\'evy process. In Section \ref{s:d-le-2} we will briefly discuss what can be said in case $d\le 2$ under somewhat stronger assumptions on the process.

Let $d\ge 3$ and let $X=(X_t,\P_x)$ be an isotropic unimodal L\'evy process in $\R^d$. This means that for each $t>0$ there is a decreasing function $p_t:(0,\infty)\to (0,\infty)$ such that
$$
\P_0(X_t\in A)=\int_A p_t(|x|)\, dx\, , \qquad A\subset \R^d \textrm{ Borel}\, .
$$
That is, transition probabilities of $X$ admit radial decreasing densities. Subordinate Brownian motions are typical examples of isotropic unimodal L\'evy processes. The characteristic exponent $\psi$ of $X$ is defined by
$$
\E_0\left[e^{i\langle x, X_t\rangle }\right]=e^{-t\psi(x)}\, , \quad x\in \R^d\, ,
$$
and is given by the L\'evy-Khintchine formula
$$
\psi(x)=a|x|^2 +\int_{\R^d\setminus\{0\}}\left(1-e^{i\langle x,y\rangle }+i\langle x, y\rangle \1_{|y|\le 1}\right)\, \nu(dy)\, ,
$$
where $a\ge 0$ and $\nu$ is the L\'evy measure having a radial decreasing density (see \cite{Wat}).
The characteristic exponent $\psi$ is a radial function: there exists $\psi_0:(0,\infty)\to (0,\infty)$ such that $\psi(x)=\psi_0(|x|)$.
In order not to overburden the notation, we will simply write $\psi(|x|)$ instead of $\psi_0(|x|)$. The same convention will be used for all radial functions on $\R^d$, e.g.~we write interchangeably $p_t(x)$ and  $p_t(|x|)$.

The main assumption we impose on the process $X$ is the following \emph{weak lower scaling condition}: There are $C_L>0$ and $\alpha>0$ such that
\begin{equation}\label{eq:lower_scaling_condition}
\psi(\lambda \xi)\geq C_L\lambda^\alpha \psi(\xi),\qquad \text{ for }\ \lambda\geq 1,\ \xi\in \R^d\, .
\end{equation}
Without loss of generality, we will assume that $C_L\le 1$.
Isotropic unimodal L\'evy processes satisfying the weak lower scaling condition were recently studied in \cite{G}.

We list several examples of such processes: (a) Brownian motion,  $\psi(x)=|x|^2$; (b) Isotropic $\beta$-stable process; $\psi(x)=|x|^{\beta}$, $\beta\in (0,2)$; (c) Independent sum of Brownian motion and isotropic $\beta$-stable process, $\psi(x)=|x|^2+|x|^{\beta}$; (d) Independent sum of two isotropic stable processes; $\psi(x)=|x|^{\alpha}+|x|^{\beta}$, $\alpha,\beta\in (0,2)$; (e) Subordinate Brownian motion via subordinator whose Laplace exponent is comparable to a regularly varying function at zero and at infinity with (not necessarily same) indices from $(0,2)$. For this class of examples see \cite{KSV8}; (f) Truncated $\beta$-stable process -- the L\'evy process with the L\'evy measure $\nu(dx)=|x|^{-d-\beta}\1_{|x|\le 1}\, dx$. For other examples see \cite[Examples 1-4]{G} and \cite{KSV8}.

By letting $\lambda \to \infty$ in \eqref{eq:lower_scaling_condition} we see that $\psi$ is unbounded. This implies that the L\'evy measure $\nu$ is infinite (unless it is zero). A consequence of this fact and transience of the process is that the Green function of $X$ is well defined and is given by
$$
G(x)=\int_0^{\infty}p_t(x)\, dt\, ,\quad x\in \R^d\, .
$$
The occupation measure of the process $X$ started at $x\in  \R^d$ is defined by $G(x,B):=\E_x\int\limits_0^\infty 1_{\{X_t\in B\}}\, dt $.
The Green function is the density of the occupation measure of the process started at the origin. Note that $G$ is radial and decreasing.

Now we are ready to state our main results which are generalizations of Theorems A and B. Recall that $A=\bigcup_{n\ge 1}\overline{B}(x_n,r_n)$ where the closed balls are pairwise disjoint and $0\notin A$. Again, we say that $A$ is avoidable (with respect to $X$) if $\P_0(T_A<\infty)<1$ where $T_A=\inf\{t>0:\, X_t\in A\}$.

\begin{Thm}\label{t:theorem-1} Let $d\ge 3$ and suppose that $X=(X_t,\P_x)$ is an isotropic unimodal L\'evy process in $\R^d$ satisfying the weak scaling condition \eqref{eq:lower_scaling_condition}.
\begin{itemize}
    \item[(a)] If $A$ is unavoidable, then
        \begin{equation}\label{e:thm-1-1}
        \sum_{n\ge 1}\frac{G(|x_n|)}{G(r_n)}=\infty\, .
        \end{equation}
    \item[(b)]Conversely, if \eqref{e:thm-1-1} and the separation condition
        \begin{equation}\label{e:thm-1-2}
        \inf\limits_{m\not= n} |x_m-x_n|^d\, \psi(|x_n|^{-1})\, G(r_n)>0
        \end{equation}
        hold, then $A$ is unavoidable.
\end{itemize}
\end{Thm}
One of the key results proved in \cite{G} states that the Green function $G(|x|)$ is comparable with $|x|^{-d}\psi(|x|^{-1})^{-1}$ (see Section 2 for more details). This means that condition \eqref{e:thm-1-1} can be written solely in terms of the given characteristic exponent $\psi$: the series in \eqref{e:thm-1-1} is divergent if and only if
$$
\sum_{n\ge 1} \frac{r_n^d \psi(r_n^{-1})}{|x_n|^d \psi(|x_n|^{-1})} =\infty\, .
$$
Moreover, we will show in Lemma \ref{l:equilibrium-potential} that $G(|x_n|)/G(r_n)$ is comparable with
$\P_0(T_{\overline{B}(x_n,r_n)}<\infty)$. Hence the series in \eqref{e:thm-1-1} is divergent if and only if
$$
\sum_{n\ge 1}\P_0(T_{\overline{B}(x_n,r_n)}<\infty)=\infty\, .
$$

\begin{Thm}\label{t:theorem-2}
Let $d\ge 3$ and suppose that $X=(X_t,\P_x)$ is an isotropic unimodal L\'evy process in $\R^d$ satisfying the weak scaling condition \eqref{eq:lower_scaling_condition}.
Suppose that the collection of balls $\{\overline{B}(x_n,r_n)\}_{n\ge 1}$ is regularly located. Then this collection is avoidable if and only if
\begin{equation}\label{e:thm-2}
\int_1^{\infty} \frac{r^{d-1}G(r)}{G(\phi(r))}\, dr <\infty\, .
\end{equation}
\end{Thm}

The proof of Theorem \ref{t:theorem-1}(a) is rather straightforward and we give it in a slightly more general form. In order to prove Theorems \ref{t:theorem-1}(b) and \ref{t:theorem-2} we modify the approach from \cite{COC} which is essentially based on estimating
the probability that starting from the origin the process $X$ exits the ball $B(0,r)$  before hitting any of the obstacles
$\overline{B}(x_n,r_n)$. An additional difficulty is caused by the fact that $X$, at the exit from $B(0,r)$, jumps out of the ball, which makes its position at the exit more uncertain than for Brownian motion which at the exit is on the boundary of $B(0,r)$. Further, since we assume the weaker separation condition \eqref{e:thm-A-2} than the one from \cite{COC}, the part using the separation condition had to be substantially modified. Here we use the scaling of the process, combined with the argument from the proof of \cite[Theorem 3]{AB}. We learned this argument from \cite{GG} who use the result of \cite[Theorem 3]{AB} in the proof of Theorem A.

An alternative way of proving Theorem \ref{t:theorem-1} is based on the concept of minimal thinness at infinity. Recall that $F\subset \R^d$ is said to be minimally thin at infinity with respect to the process $X$ if $\P_0(T_F<\infty)<1$. Thus, $A=\bigcup_{n\ge 1}\overline{B}(x_n, r_n)$ is avoidable if and only if it is minimally thin at infinity. A Wiener-type criterion (see \cite{AE}) for minimal thinness at infinity with respect to Brownian motion is well known and can be described as follows: For $j\in \Z$ consider the closed cube centered at $0$, with sidelength $3^j$, and sides  parallel to coordinate axes. Divide this cube into $3^d$ subcubes of sidelength $3^{j-1}$ and discard the central cube. The enumeration  $\{Q_m\}_{m\ge 1}$ of all such cubes is a Whitney decomposition of
$\R^d\setminus \{0\}$. Let $\mathrm{Cap}$ denote the usual Newtonian capacity. Suppose that $F$ is closed and does not contain the origin. Then $F$ is minimally thin at infinity (with respect to Brownian motion) if and only if
$$
\sum_{m\ge 1}\mathrm{diam}(Q_m)^{2-d}\, \mathrm{Cap}(F\cap Q_m) <\infty\, .
$$
This criterion is used in \cite{GG} in order to prove (the converse part of) Theorem A (see \cite[Theorem 6]{GG}). It is rather straightforward, although technically demanding, to extend the latter approach to the class of isotropic unimodal L\'evy processes satisfying the weak lower scaling condition \eqref{eq:lower_scaling_condition}. The first step in this approach would be an analog of the Wiener-type criterion for minimal thinness at infinity with respect to the L\'evy process $X$ in terms of the corresponding capacity. Here we need an additional assumption on the process, namely that it is a subordinate Brownian motion. Analytically this means that the characteristic exponent $\psi$ is of the form $\psi(x)=f(|x|^2)$ where $f$ is a Bernstein function.
Below $\cp$ denotes the capacity with respect to the process $X$.

\begin{Thm}\label{t:aikawa-criterion}
Let $d\ge 3$ and suppose that $X=(X_t,\P_x)$ is a subordinate Brownian motion in $\R^d$ satisfying the weak lower scaling condition \eqref{eq:lower_scaling_condition}. Let $F\subset \R^d\setminus \{0\}$ be a closed set. Then $F$ is minimally thin at infinity with respect to $X$ if and only if
\begin{equation}\label{e:aikawa-criterion}
\sum_{m\ge 1}G(\mathrm{diam}(Q_m))\, \mathrm{Cap}(F\cap Q_m) <\infty\, .
\end{equation}
\end{Thm}
The key step in proving this theorem is to show quasiadditivity of the capacity with respect to the Whitney decomposition $\{Q_m\}_{m\ge 1}$. This can be achieved by modifying the method employed in \cite{Aik} and \cite{AE}. Since the full proof would be quite long, we only outline the main steps and changes in the Appendix.

Theorem \ref{t:aikawa-criterion} can be used to study avoidability of a collection of balls with random centers given by the Poisson point process. In case of Brownian motion this question was recently studied in \cite{CDOC}. By adopting their method we are able to 
extend the results to subordinate Brownian motions satisfying \eqref{eq:lower_scaling_condition}. We state the assumptions and
the result in Section \ref{s:poisson}. Again, only  necessary changes in the proof are outlined.

In Section \ref{s:d-le-2} we look at the case $d\le 2$. Here we need somewhat stronger assumptions on the process. We will assume that $X$ is a transient subordinate Brownian motion via the subordinator whose Laplace exponent is a complete Bernstein function which satisfies the lower and the upper scaling conditions at zero and infinity. These conditions, called {\bf (H1)} and {\bf (H2)}, were used in \cite{KSV8} to study some aspects of the potential theory. With these conditions in place, Theorems \ref{t:theorem-1} and \ref{t:theorem-2} are still valid, and we indicate modifications of proofs. The reason stronger assumptions are needed can be explained as follows: The weak lower scaling condition \eqref{eq:lower_scaling_condition} gives the lower bound on the quotient $\psi(\lambda \xi)/\psi(\xi)$. A corresponding upper bound $\psi(\lambda \xi) \le 4\lambda^2 \psi(\xi)$ - see \eqref{eq:square} - is always true,  but with $\lambda^2$ instead of $\lambda^{\alpha}$. This upper bound is sufficient for our purposes in case $d\ge 3$. This is not so when $d\le 2$, which makes it necessary to impose an appropriate upper bound.

Finally, we make a comment on the related, but more difficult, problem of avoiding a collection of balls that are contained in the unit ball. In case of Brownian motion, this has been studied in \cite{Ake, Don, GG, HN, OCS}. 
For L\'evy processes
such problems make no sense since with positive probability the process jumps out of the unit ball, thus making every collection avoidable.
A natural jump processes for the avoidability problem in bounded open set $D$ is a censored $\alpha$-stable process with $\alpha\in (1,2)$. This is a process obtained from isotropic stable process by suppressing jumps landing outside of $D$ and continuing at the place where the suppressed jump has occurred. Such process is transient and approaches the boundary at its lifetime. Avoidability of collection of balls in bounded $C^{1,1}$ open set is studied in \cite{AM}. In a more general context of balayage spaces the problem of smallness of unavoidable sets has been addressed in \cite{HN13}. 

We finish this introduction with a few words on notation. Throughout the paper we use a number of constants. Constants whose values are important will be denoted by uppercase letter $C_1, C_2,\dots$, or will have subscripts to remind the reader where they come from, e.g., the already introduced $C_L$ in \eqref{eq:lower_scaling_condition}. Unimportant constants will be denoted by $c_1, c_2, \dots $. Throughout the paper we use the notation $f(r) \asymp g(r)$ as  $r \to a$ to denote that $f(r)/g(r)$ stays between two positive constants as $r \to a$.

\section{Preliminary Results}\label{s:preliminaries}
Let $d\ge 3$ and let $X=(X_t,\P_x)$ be an isotropic unimodal L\'evy process in $\R^d$ with the characteristic exponent $\psi$ satisfying the weak lower scaling condition \eqref{eq:lower_scaling_condition}.

The characteristic exponent $\psi$ is a radial function. Following \cite{G}, we introduce an increasing modification  $\psi^{\star}:[0,\infty)\to [0,\infty)$ of $\psi$ by
$$
\psi^\star(r)=\sup_{s\in [0,r]}\psi(s),\quad r> 0\, .
$$
By \cite[Lemma 1]{G}, $\psi^{\star}$ satisfies the following weak lower scaling condition:
\begin{equation}\label{e:lsc-2}
\psi^{\star}(\lambda r)\ge C_L^2 \lambda^{\alpha}\psi^{\star}(r)\, ,\qquad  \lambda\geq 1,\ r>0\, .
\end{equation}
By \cite[Proposition 1]{G}, for any $r>0$,
\begin{equation}\label{eq:comp}
	\psi(r)\leq \psi^\star(r)\leq 12\psi(r)\, .
\end{equation}
Further, by \cite[Lemma 1]{G}
\begin{eqnarray}
\psi^\star(\lambda r)&\leq & 2(1+\lambda ^2)\psi^\star(r)\, ,\quad \lambda, r>0\, ,\label{eq:square}\\
\psi^\star(\lambda r) &\geq &\frac12 \frac{\lambda^2}{1+\lambda^2}\psi^\star(r)\, ,\quad \lambda, r>0\, .\label{eq:square-left}
\end{eqnarray}	
When $\lambda\ge 1$, \eqref{eq:square} reads as weak upper scaling condition on $\psi^{\star}$: $\psi^\star(\lambda r)\le 4\lambda^2 \psi^{\star}(r)$, $r>0$. Note the power $2$ on $\lambda$ which will be sufficient when $d\ge 3$.

We will need the following two estimates.
\begin{Lem}\label{lem:green_est}
There exists a constant $C_G=C_G(d, C_L, \alpha)\geq 1$ such that
	\begin{equation}\label{e:green-est}
		\frac{C_G^{-1}}{|x|^d\psi^\star (|x|^{-1})}\leq G(x)\leq \frac{C_G}{|x|^d\psi^\star (|x|^{-1})}\, ,\qquad x\in \R^d\setminus\{0\}\,.
	\end{equation}
In particular,
\begin{equation}\label{eq:green_scaling}
	\frac{1}{4C_G^2}\lambda^{2-d}G(x)\leq G(\lambda x)\leq \frac{C_G^2}{C_L}\lambda^{\alpha-d}G(x)\, ,\qquad \lambda\in (0,1]\, ,\  x\in \R^d\setminus\{0\}\, .
\end{equation}
\end{Lem}
\proof
Inequality \eqref{e:green-est} is proved in \cite[Theorem 3]{G}. By (\ref{eq:lower_scaling_condition}) and (\ref{eq:square}),
$$
	C_L\lambda^{-\alpha}\psi^\star(|x|^{-1})\leq \psi^\star(\lambda^{-1}|x|^{-1})\leq 4\lambda^{-2}\psi^\star(|x|^{-1})\,,
$$
which implies inequality (\ref{eq:green_scaling}).
\qed

Note that as the special case of the right-hand side inequality in \eqref{eq:green_scaling} (with $\lambda=1/2$) we get the following doubling condition for the Green function: There is a constant $c>0$ such that
\begin{equation}\label{e:doubling-G}
G(2x)\ge c\ G(x)\, , \quad x\in \R^d\setminus\{0\}\, .
\end{equation}

For an open set $D\subset \R^d$ we define the \emph{first exit time} from $D$ by $\tau_D=\inf\{t>0\colon X_t\not\in D\}$\,.

\begin{Lem}{\cite[Corollary 2]{G}}\label{l:long-jump}
	There is a constant $C_J=C_J(d)>0$ such that for any $x_0\in \R^d$ and $r,s>0$, $s\geq 2r$,
	$$
		\sup_{x\in B(x_0,r)}\P_x(X_{\tau_{B(x_0,r)}}\not\in B(x_0,s))\leq C_J \left(\frac{r}{s}\right)^\alpha\,.
	$$
\end{Lem}

\begin{Def}
A function $u\colon \R^d\rightarrow \R$ is {\it harmonic} in an open set $D\subset \R^d$ if for any open set $B\subset \overline{B}\subset D$ and $x\in B$ the following holds
	\begin{equation}\label{eq:maxpr}
		\E_x[|u(X_{\tau_B})|; \tau_B<\infty]<\infty\qquad\text{ and }\qquad u(x)= \E_x[u(X_{\tau_B}); \tau_B<\infty]\,.
	\end{equation}
We say that $u$ is \emph{regular harmonic} in $D$ if (\ref{eq:maxpr}) holds with $B=D$\,.
\end{Def}
For basic properties of harmonic functions in this context see for example \cite{SV04}.

\begin{Lem}[Maximum Principle]\label{lem:maxpr}
	Let $D\subset \R^d$ be an open set and let $u,v\colon \R^d\rightarrow \R$ be regular harmonic functions in $D$ satisfying
	$$
		u(x)\geq v(x),\qquad x\in D^c\, .
	$$
	Then $u(x)\geq v(x)$ for all $x\in \R^d$.
\end{Lem}
\proof
	Since $w=u-v$ is regular harmonic in $D$ and $w(y)\geq 0$ for $y\in D^c$ we get
	$$
		w(x)=\E_x[w(X_{\tau_D});\tau_D<\infty]\geq 0,\quad x\in D\, .
	$$
\qed

The \emph{first hitting time} of a closed set $F\subset \R^d$ is defined by $T_F=\inf\{t>0\colon X_t\in F\}$.

\begin{Lem}\label{l:equilibrium-potential}
	There is a constant $C_E\in (0,1]$ so that  for all $x_0\in \R^d$, $r>0$ and $x\in B(x_0,r)^c$,
    $$
        C_E \frac{G(x-x_0)}{G(r)}\leq \P_x(T_{\overline{B}(x_0,r)}<\infty)\leq  \frac{G(x-x_0)}{G(r)}\,.
	$$
\end{Lem}
\proof
Without loss of generality we assume that $x_0=0$. Set $B=\overline{B}(0,r)$. Since $G$ and $x\mapsto \P_x(T_B<\infty)$ are regular harmonic in  $B^c$ and $\P_x(T_B<\infty)=1\leq \frac{G(x)}{G(r)}$ for all
$x\in B$, by Lemma \ref{lem:maxpr}  we get the upper bound.

Let $x\in B^c$. By using the strong Markov property we obtain
	\begin{align*}
		G(x,B)&=\E_x\int\limits_0^\infty 1_{\{X_t\in B\}}\,dt=\E_x\left[\int_{T_B}^\infty 1_{\{X_t\in B\}}\,dt;T_B<\infty\right]
        \\&=\E_x[G(X_{T_B},B);T_B<\infty]\,.
	\end{align*}
 Since $G(z,B)\leq \int\limits_{B(0,2r)}G(y)\, dy$ for $z\in B$, it follows that
	$$
		\P_x(T_B<\infty)\geq \frac{\int\limits_B G(y-x)\,dy}{\int\limits_{B(0,2r)}G(y)\,dy}\,.
	$$
By Lemma \ref{lem:green_est}, \eqref{eq:comp}, lower scaling condition \eqref{eq:lower_scaling_condition} and  \eqref{eq:square} it follows that
	$$
		\int\limits_{B(0,2r)}G(y)\,dy\leq c_1 \int\limits_0^{2r}\frac{ds}{s\psi^\star(s^{-1})}\leq c_2\psi^\star((2r)^{-1})^{-1} \int\limits_0^{2r}\frac{ds}{s\left(\frac{2r}{s}\right)^\alpha}\leq c_3 \psi^\star(r^{-1})^{-1}\, .
	$$
Since $|y-x|\leq |y|+|x|\le r+|x| \le 2|x|$ for all $y\in B$, the doubling property of the Green function implies
	$$
		\int\limits_B G(y-x)\,dy\geq G(2x)|B|\geq c_4 G(x)r^d\,.
	$$
    Hence, from the last three displayed equations and Lemma \ref{lem:green_est} we deduce,
	\begin{equation}\label{eq:exit1}
        \P_x(T_B<\infty)\geq \frac{c_4}{c_3}\frac{G(x)r^d}{\psi^\star(r^{-1})^{-1}}\geq c_5 \frac{G(x)}{G(r)},\quad x\in B(0,r)^c\, .
	\end{equation}
\qed

We now discuss scaling of the process $X$. First note that without loss of generality we may assume that $\psi^{\star}(1)=\psi(1)$. For $a>0$ let
$$
\psi^a(\xi):=\frac{\psi(a\xi)}{\psi^{\star}(a)}\, ,\quad \xi\in \R^d\, .
$$
Then clearly,
$$
\psi^{a,\star}(r):=\sup_{s\in [0,r]}\psi^a(s)=\frac{\psi^{\star}(ar)}{\psi^{\star}(a)}\, ,\quad r>0\, .
$$
Since $\psi^a$ is a continuous, negative definite function, it is a characteristic exponent of a L\'evy process $X^a=(X^a_t)_{t\ge 0}$. It is easy to see that $X^a$ has the same law as the scaled process $(aX_{t/\psi^{\star}(a)})_{t\ge 0}$. Indeed
$$
\E_0\left[e^{i\langle x, aX_{t/\psi^{\star}(a)}\rangle}\right]=e^{-\frac{t}{\psi^{\star}(a)}\, \psi(ax)}=e^{-t\psi^a(x)}\, ,\quad x\in \R^d\, .
$$
Further, for $\lambda\ge 1$ and $\xi \in \R^d$,
$$
\psi^a(\lambda \xi)=\frac{\psi(\lambda a \xi)}{\psi^{\star}(a)}\ge C_L \lambda^{\alpha}\frac{\psi(a\xi)}{\psi^{\star}(a)} =C_L\lambda^{\alpha}\psi^a(\xi)\, .
$$
Hence,  $X^a$ is also isotropic unimodal and satisfies the same weak lower scaling condition \eqref{eq:lower_scaling_condition}. From the equality in law of $X^a$ and $(aX_{t/\psi^{\star}(a)})_{t\ge 0}$, it immediately follows that for any Borel set $B$ it holds that
\begin{equation}\label{e:equilibrium-scaling}
\P_x(T_B<\infty)=\P_{ax}(T^a_{aB}<\infty)\, ,\quad x\in \R^d\, ,
\end{equation}
where $aB=\{ay:\, y\in B\}$ and $T^a_{aB}=\inf\{t>0:X^a_t\in aB\}$ is the hitting time of $B$ by $X^a$\,.
\begin{Lem}\label{lem:scaling}
Let $a>0$ and let $X=(X^a_t)_{t\geq 0}$ be a L\'evy process with the L\'evy exponent $\psi^a$. The Green function $G^a$ of $X^a$ is given by
		\begin{equation}\label{eq:green_scale2}
			G^a(x)=a^{-d}\psi^\star(a)G(a^{-1}x)\, ,\quad x\in \R^d\, .
		\end{equation}
In particular, if the weak scaling condition \eqref{eq:lower_scaling_condition} holds, then for every $\rho>0$ there is a constant $C_1=C_1(d,C_L,\alpha,\rho)>0$ such that
		$$
			G^a(x)\leq C_1\left(\,|x|^{2-d}\vee |x|^{\alpha-d}\right),\qquad |x|\leq \rho\, .
		$$
\end{Lem}
\proof Let $B\subset \R^d$ be a Borel set. By equality in law of $X^a$ and $(aX_{t/\psi^{\star}(a)})_{t\ge 0}$ and a change of variables, we see that
	\begin{align*}
		\int\limits_0^\infty \P_0(X^a_t\in B)\,dt&=\int\limits_0^\infty \P_0(aX_{t/\psi^{\star}(a)}\in B)\,dt
		=\psi^\star(a)\int\limits_0^\infty \P_0(X_t\in a^{-1}B)\,dt\\
		&=\psi^\star(a)\int\limits_{a^{-1}B} G(y)\,dy
        =\int\limits_{B} a^{-d}\psi^\star(a)G(a^{-1}y)\,dy.
	\end{align*}
This shows that the Green function of $X^a$ exists and is given by \eqref{eq:green_scale2}.

Assume now that \eqref{eq:lower_scaling_condition} holds.
By Lemma \ref{lem:green_est},
	$$
		G^a(x)\leq C_G|x|^{-d}\frac{\psi^\star(a)}{\psi^\star(a|x|^{-1})}\, .
	$$
Let $N\in \N\cup\{0\}$ be such that $\rho^{-1}\geq 2^{-N}$, i.e.~$N\geq \log_2{\rho}$\,.
Then $|x|^{-\frac{1}{N}}\geq \rho^{-\frac{1}{N}}\geq 1/2$. Suppose that $|x|^{-1/N}\ge 1$ (i.e.~$|x|\le 1$). Then by \eqref{eq:lower_scaling_condition},
	$$
    \frac{\psi^\star(a)}{\psi^\star(a|x|^{-1})} =\frac{\psi^\star(a)}{\psi^\star(a|x|^{-\frac{1}{N}})}\frac{\psi^\star(a|x|^{-\frac{1}{N}})}{\psi^\star(a|x|^{-\frac{2}{N}})} \cdots\frac{\psi^\star(a|x|^{-\frac{N-1}{N}})}{\psi^\star(a|x|^{-1})}\leq C_L^{-N}|x|^{\alpha}\, .
	$$
If $1/2\le |x|^{-1/N}\le 1$ (i.e.~$1\le |x|\le 2^N$) we use \eqref{eq:square-left} to conclude that for every $b>0$
$$
\psi^{\star}(|x|^{-1/N}b)\ge \frac12 \, \frac{|x|^{-2/N}}{1+|x|^{-2/N}}\, \psi^{\star}(b) \ge \frac14 |x|^{-2/N} \psi^{\star}(b)\, .
$$
Now the same computation as above gives that
$$
\frac{\psi^\star(a)}{\psi^\star(a|x|^{-1})}\le 4^N |x|^2\, .
$$
	Finally,
	$$
		G^a(a)\leq (C_L^{-N}+4^N) C_G \left(\,|x|^{2-d}\vee |x|^{\alpha-d}\right),\qquad |x|\leq \rho\, .
	$$
\qed

Let $\cp$ denote the capacity with respect to $X$ and $\cp^a$ the capacity with respect to $X^a$. It is shown in \cite[Proposition 3]{G} that there exists a constant $C_2=C_2(d)>1$ such that for any $r>0$,
\begin{equation}\label{e:capacity-of-ball}
C_2^{-1}\psi^{\star}(r^{-1})r^d \le \cp(\overline{B}(0,r))\le C_2 \psi^{\star}(r^{-1})r^d\, .
\end{equation}

\begin{Lem}\label{l:capacity-estimate}
There exists a constant $C_3=C_3(d,C_L,\alpha)>1$ such that for every $a>0$ and every $r>0$
$$
C_3^{-1}\frac{1}{G^a(r)}\le \cp^a(\overline{B}(0,r)) \le C_3 \frac{1}{G^a(r)}\, .
$$
\end{Lem}
\proof Since $X^a$ satisfies the same assumptions as $X$ with same constants, it suffices to prove the statement for $a=1$. By
\cite[Proposition 5.55]{SV}, there exists a constant $c_1=c_1(d)>1$ such that for every $r>0$
\begin{equation}\label{e:cap-est-1}
c_1^{-1}\frac{r^d}{\int_{B(0,r)}G(y)\, dy}\le \cp(\overline{B}(0,r)) \le c_1 \frac{r^d}{\int_{B(0,r)}G(y)\, dy}\, .
\end{equation}
By using \eqref{e:green-est} and polar coordinates we see that $\int_{B(0,r)}G(y)\, dy \asymp \int_0^r s^{-1}\psi^{\star}(s^{-1})^{-1}\, ds$ with constants depending on $d$, $C_L$ and $\alpha$. The weak lower scaling condition implies the upper bound $\int_0^r s^{-1}\psi^{\star}(s^{-1})^{-1}\, ds \le c_2 \psi^{\star}(r^{-1})^{-1}$ with $c_c=c_2(C_L, \alpha)$, while \eqref{eq:square} implies the lower bound $\int_0^r s^{-1}\psi^{\star}(s^{-1})^{-1}\, ds \ge c_3 \psi^{\star}(r^{-1})^{-1}$. By inserting in \eqref{e:cap-est-1} and using \eqref{e:green-est} again, we get that
$$
\cp(\overline{B}(0,r))\asymp r^d \psi^{\star}(r^{-1})\asymp \frac{1}{G(r)}\, .
$$
\qed

\begin{Lem}\label{l:capacity-scaling}
For any bounded Borel set $B\subset \R^d$ and for every $a>0$ it holds that
$$
\cp^a (aB)=a^d\psi^{\star}(a)^{-1}\cp(B)\, .
$$
\end{Lem}
\proof Let $\mu$ denote the equilibrium measure of the set $B$, i.e., $\P_x(T_B<\infty)=G\mu(x)$. Then $\cp(B)=\mu(B)$ (see, e.g.\cite[VI.4]{BG}). Let $\mu^a$ denote the equilibrium measure of $aB$ with respect to the process $X^a$, and define the measure $\wt{\mu}$ by $\wt{\mu}(A):=a^{-d}\psi^{\star}(a)\mu^a(aA)$. By using \eqref{eq:green_scale2} in the second, and change of variables in the third line, we have
\begin{eqnarray*}
\P_{ax}(T^a_{aB}<\infty)&=&G^a \mu^a(ax)=\int G^a(y-ax)\, \mu^a(dy)\\
&=&a^{-d}\psi^{\star}(a)\int G(a^{-1}y-x)\, \mu^a(dy)\\
&=&\int G(y-x)\, \wt{\mu}(dy)=G\wt{\mu}(x)\, .
\end{eqnarray*}
On the other hand, $\P_{ax}(T^a_{aB}<\infty)=\P_x(T_B<\infty)=G\mu(x)$. By the uniqueness principle (see \cite[VI.1 Proposition 1.15]{BG}) we conclude that $\wt{\mu}=\mu$. Hence,
$$
\cp(B)=\mu(B)=\wt{\mu}(B)=a^{-d}\psi^{\star}(a)\mu^a(aB)=a^{-d}\psi^{\star}(a)\cp^a(aB)\, ,
$$
proving the claim. \qed

\section{Proof of Theorem \ref{t:theorem-1}}\label{s:proof-1}
We start with a general simple sufficient condition for avoidability of collection of balls from which the first part of Theorem \ref{t:theorem-1} will immediately follow. Let $\{\overline{B}(x_n,r_n)\}_{n\ge 1}$ be a family of disjoint closed balls in $\R^d$ and $A:=\bigcup_{n\ge 1}\overline{B}(x_n,r_n)$.

\begin{Prop}\label{p:sufficient-cond}
Let $d\ge 3$ and suppose that $X=(X_t,\P_x)$ is an isotropic unimodal L\'evy process in $\R^d$ satisfying the weak scaling condition \eqref{eq:lower_scaling_condition}. If
    \begin{equation}\label{e:sufficient-condition}
    \sum_{n\ge 1}\P_0(T_{\overline{B}(x_n,r_n)}<\infty)<\infty\, ,
    \end{equation}
then $A$ is avoidable.
\end{Prop}
\proof Recall first that $0\notin A$. To simplify notation, let $B_n=\overline{B}(x_n,r_n)$, $n\ge 1$. Choose $N\in \N$ such that $\sum\limits_{n>N}^{\infty}\P_0(T_{B_n}<\infty)<1/2$, and let $G:=\bigcup\limits_{n>N} B_n$. Since $\{T_G<\infty\}\subset \bigcup_{n>N}\{T_{B_n}<\infty\}$, we have that
$$
\P_0(T_G<\infty)\le \sum\limits_{n>N}\P_0(T_{B_n}<\infty)<\frac12\, .
$$
Let $u(x):=\P_x(T_G<\infty)$, $x\in \R^d$. The function $u$ is regular harmonic in $G^c$, and by the above, $u(0)<1/2$. We claim that the set $\{x\in G^c:\, u(x)<1/2\}$ is unbounded. If not, $D:=\{x\in G^c\colon u(x)<1/2\}$ is bounded, hence $\P_0(\tau_D<\infty)=1$.
Since $u\equiv1$ on $G$, $u(x)\geq \frac{1}{2}$ on $D^c$. Regular harmonicity of $u$ in $G^c$ implies that $u(0)=\E_0[u(X_{\tau_D})]\ge 1/2$. Contradiction!

Now define $F:=\bigcup\limits_{n=1}^N B_n$.
Since $F$ is bounded, the right-hand side inequality in Lemma \ref{l:equilibrium-potential} implies that $\lim\limits_{|x|\to\infty}\P_x(T_F<\infty)=0$\,.
Hence, there exists a point $y\in G^c$ such that $\P_y(T_F<\infty)<1/2$ and $\P_y(T_G<\infty)<1/2$. Since $A=F\cup G$, we have that $$
\P_y(T_A<\infty)\le \P_y(T_F<\infty)+\P_y(T_G<\infty)<\frac12 +\frac12=1\, .
$$
Let $v(x):=\P_x(T_A<\infty)$. Then $v$ is bounded in $\R^d$ and regular harmonic in $A^c$. It is proved in \cite[Theorem 2]{G} that bounded harmonic functions are (H\"older) continuous. Define  $U:=\{v<1\}$. Then $U$ is nonempty (since  $y\in U)$ and open.

In the sequel we distinguish whether $X$ has a jump component or not. If it does not, then $X$ is Brownian motion and the (classical) maximum principle implies that $U=A^c$ (recall that $A^c$ is connected). In case there is a jump component, i.e.~the L\'evy measure is non-trivial, we will prove that $\P_0(T_U<T_A)>0$. With this result we proceed in the following way. Since $(v(X_t^{T_A}))_{t\ge 0}$ is a $\P_0$-martingale, by the optional stopping theorem we have for all $t>0$
\begin{eqnarray*}
v(0)&=&\E_0[v(X_{T_U\wedge T_A\wedge t})]=\E_0[v(X_{T_U\wedge t}), T_U<T_A]+\E_0[v(X_{T_A\wedge t}), T_A\le T_U]\\
&\le &\E_0[v(X_{T_U\wedge t}), T_U<T_A] +\P_0(T_A\le T_U)\, .
\end{eqnarray*}
By letting $t\to \infty$, the first term above converges to $\E_0[v(X_{T_U}), T_U<T_A]<\P_0(T_U<T_A)$. Thus $v(0)<1$, proving that $A$ is avoidable.

It remains to show that $\P_0(T_U<T_A)>0$. Let $j$ denote the radially decreasing density of the L\'evy measure $\nu$: $\nu(dx)=j(x)\, dx$. Suppose first that $\nu$ has unbounded support. Let $\epsilon>0$ be such that $B(0,\epsilon)\cap B(y,\epsilon)=\emptyset$, $B(y,\epsilon)\subset U$, and $A$ does not intersect either of the balls. Then by the Ikeda-Watanabe formula,
$$
\P_0(T_U<T_A)\ge \P_0(X_{\tau_{B(0,\epsilon)}}\in B(y,\epsilon))= \int_{B(y,\epsilon)}\int_{B(0,\epsilon)}G_{B(0,\epsilon)}(x,dy)j(z-y)\, dz >0\, .
$$
The general case is slightly more complicated. Suppose that $j(x)>0$ for $|x|\le R$, $R>0$. Let $0=x_0,x_1,\dots, x_n=y$ be a sequence of points in $A^c$ such that for $\epsilon \in (0,R/8)$ small enough, the balls $B(x_j, \epsilon)$ are pairwise disjoint and contained in $A^c$, and moreover $|x_k-x_{k-1}|<R/4$ for all $k=1,2, \dots, n$. We first show that
$$
c:=\inf_{x\in B(x_{k-1}, \epsilon/2)}\P_x(X_{\tau_{B(x_{k-1},\epsilon)}}\in B(x_k,\epsilon/2))>0\, .
$$
Indeed, by the Ikeda-Watanabe formula, for every $x\in B(x_{k-1}, \epsilon/2)$ we have
\begin{eqnarray*}
\P_x(X_{\tau_{B(x_{k-1},\epsilon)}}\in B(x_k,\epsilon/2))& = &\int_{B(x_k, \epsilon/2)} \int_{B(x_{k-1}, \epsilon)}G_{B(x_{k-1}, \epsilon)}(x,z) j(z-w)\, dz\, dw\\
&\ge & j(R/2)|B(x_k, \epsilon/2)|\, \int_{B(x_{k-1}, \epsilon)}G_{B(x_{k-1}, \epsilon)}(x,z)\, dz\\
& =& j(R/2)|B(x_k, \epsilon/2)|\, \E_x \tau_{B(x_{k-1}, \epsilon)}\, .
\end{eqnarray*}
Since $\inf_{x\in B(x_{k-1},\epsilon/2)} \E_x \tau_{B(x_{k-1}, \epsilon)}>0$, we get that $c>0$. Now, $\P_0(T_U<T_A)$ is certainly larger then the probability of successively exiting the balls $B(x_{k-1},\epsilon)$ by jumping into $B(x_k,\epsilon/2)$, $k=1,2,\dots, n$. The latter probability is bounded from below by $c^n$. Thus $\P_0(T_U<T_A)\ge c^n>0$.
\qed

\emph{Proof of Theorem \ref{t:theorem-1}(a)}: The statement follows directly from Proposition \ref{p:sufficient-cond} and Lemma \ref{l:equilibrium-potential} which states that $\P_0(T_{\overline{B}(x_n,r_n)}<\infty) \asymp G(|x_n|)/G(r_n)$. \qed

\begin{Lem}\label{lem:recurrence}
Let $B\in \B(\R^d)$ and let $V_1,V_2,W\subset \R^d$ be open sets such that $\overline{V_1}\subset W\subset \overline{W}\subset V_2$. Then
	$$
        \P_0(\tau_{V_2}<T_B)\leq \P_0(\tau_{V_1}<T_B)\left(t\sup_{y\in W\setminus V_1}\P_y(\tau_{V_2}<T_B)+1-t\right),
	$$
	where
	$$
        t:=\P_0\left(X_{\tau_{V_1}}\in W\, \big|\, \tau_{V_1}<T_B\right)\, .
	$$
\end{Lem}
\proof
Suppose that $0\in V_1$ which is the case of interest in the forthcoming results.  It is easily checked that the statement holds true also for $0\not\in V_1$.
Separating the event in the probability on the left-hand side and using the strong Markov property it follows that
	\begin{eqnarray*}
	\lefteqn{\P_0(\tau_{V_2}<T_B)=\P_0(\tau_{V_2}<T_B, X_{\tau_{V_1}}\in W)+\P_0(\tau_{V_2}< T_B, X_{\tau_{V_1}}\not\in W)}\\
	&\leq &\E_0[\P_{X_{\tau_{V_1}}}(\tau_{V_2}<T_B);\tau_{V_1}<T_B, X_{\tau_{V_1}}\in W]+\P_0(\tau_{V_1}<T_B, X_{\tau_{V_1}}\not\in W)\\
		&\leq & \sup_{y\in W\setminus V_1}\P_y(\tau_{V_2}<T_B) \P_0(\tau_{V_1}<T_B, X_{\tau_{V_1}}\in W)
    +\P_0(\tau_{V_1}<T_B, X_{\tau_{V_1}}\not\in W)\, .\quad \qed
	\end{eqnarray*}

\proof[Proof of Theorem \ref{t:theorem-1}(b)]
Assume that the separation condition \eqref{e:thm-1-2} holds and that $A$ is avoidable. We will prove that
$$
\sum_{n\ge 1}\frac{G(|x_n|)}{G(r_n)}<\infty\, .
$$
First note that without loss of generality we may assume that $r_k\le |x_k|/2$ for all $k\ge 1$.
Indeed, if we put $r_k'=r_k\wedge |x_k|/2$, and the collection $\{\overline{B}(x_k,r_k)\}_{k\ge 1}$ is avoidable, then clearly the collection $\{\overline{B}(x_k,r_k')\}_{k\ge 1}$ is avoidable as well, and condition \eqref{e:thm-1-2} is true with $r_k'$. Suppose that we have proved that $\sum_{n\ge 1}G(|x_n|)/G(r_n')<\infty$. Then $r_n'=r_n$ for all $n$ sufficiently large. Otherwise $r_n'=|x_n|/2$ for infinitely many $n$, hence $G(r_n')\asymp G(|x_n|)$ for infinitely many $n$, and the series will diverge.

Let $c_0\in (0,1)$ be such that
\begin{equation}\label{eq:cond_variant}
    \inf\limits_{j\not= k} |x_j-x_k|^d\psi^\star(|x_k|^{-1})\, G(r_k)>c_0\, .
\end{equation}
Pick $\rho \geq 2$ large enough so that
\begin{equation}\label{e:choice-rho}
	4C_G^2\left(\frac{2\rho}{\rho^2-1}\right)^{d-2}\leq\frac{C_E}{2}\,.
\end{equation}
Set
\begin{equation}\label{eq:choice_nu_delta}
	\nu:=\left(\frac{2^{d}}{c_0 C_L}\right)^{1/\alpha}\geq 1 \quad\text{ and }\quad
    \delta:=\left(\frac{c_0C_L}{4^{d}C_1}\right)^{1/\alpha}\frac{1}{\rho}\leq 1,
\end{equation}
where $C_1=C_1(\nu \rho)\geq 1$ is the constant from Lemma \ref{lem:scaling}.
This choice of $\delta$ implies that for all $r\le \nu \rho$,
\begin{equation}\label{eq:green_scale}
	G(\delta r)\leq C_1(\delta r)^{\alpha-d}\leq 2^{-d}(\delta r)^{-d}\,.
\end{equation}
Note that $G(\delta r_n)\asymp G(r_n)$ (with a constant independent of $n$). This, together with
$\delta\leq 1$, implies that it suffices to prove the theorem for the balls $\{\overline{B}(x_n,\delta r_n)\}_{n\geq 1}$. Therefore, in the sequel we set $\tilde{r}_n:=\delta  r_n$ and consider balls $\{\overline{B}(x_n,\tilde{r}_n)\}_{n\geq 1}$\,.

Let $\theta:=\P_0(T_A=\infty)>0$ by the assumption and suppose that the series in \eqref{e:thm-1-1} is divergent.  Let $B_n:=B(0, \rho^n)$ and $p_n:=\P_0(\tau_{B_n}<T_A)$, $n\ge 1$. We are going to prove that $\lim\limits_{n\to \infty}p_n=0$, contradicting $\P_0(T_A=\infty)>0$.

We now chose $m_0\in \N$ large enough so that
\begin{equation}\label{e:choice-m0}
\frac{C_J}{\rho^{(m_0-2)\alpha}}\le \frac{\theta}{2}\ \quad \text{ and }\quad C_L\rho^{\alpha m_0}>12\,.
\end{equation}

Set $I_k:=\{n\in \N:\,  \rho^{k-1}<|x_n|\le  \rho^k\}$, $k\in \N$.
Since we have assumed that the series in \eqref{e:thm-1-1} diverges, there exists $\ell\in \{0,1,\dots, m_0-1\}$ such that
\begin{equation}\label{e:smaller-series}
\sum_{j=0}^{\infty}\sum_{n\in I_{\ell+m_0 j}} \frac{G(|x_n|)}{G(\tilde{r}_n)}=+\infty\, .
\end{equation}
Set
$$
m_j:=\ell+m_0 j\, \qquad q_j:=\sup_{y\in B_{m_{j+1}-2}\, \setminus B_{m_j}}\P_y (\tau_{B_{m_{j+1}}}<T_A)\, ,
$$
and
$$
 \quad t_j:=\P_0\left(X_{\tau_{B_{m_j}}}\in B_{m_{j+1}-2}\,  \big|\, \tau_{B_{m_j}}<T_A\right)\, .
$$
By Lemma \ref{lem:recurrence} (with $V_1=B_{m_j}$, $W=B_{m_{j+1}-2}$, $V_2=B_{m_{j+1}}$),
\begin{equation}\label{e:recursion}
p_{m_{j+1}}\le p_{m_j}(q_j t_j+1-t_j)\, ,\quad j\in \N\, .
\end{equation}
We will prove that
\begin{equation}\label{e:qj-diverge}
\sum_{j=1}^{\infty}(1-q_j)=+\infty\, .
\end{equation}
Suppose that \eqref{e:qj-diverge} is true. Then
$$
t_j\ge \frac{p_{m_j}-\P_0(X_{\tau_{B_{m_j}}}\notin B_{m_{j+1}-2})}{p_{m_j}}\ge \frac{\theta -\frac{\theta}{2}}{p_{m_j}}\ge \frac{\theta}{2}
$$
by Lemma \ref{l:long-jump} and condition \eqref{e:choice-m0}. Hence,
$$
\sum_{j=1}^{\infty}\big(1-(q_j t_j+1-t_j)\big)=\sum_{j=1}^{\infty}(1-q_j)t_j \ge \frac{\theta}{2}\sum_{j=1}^{\infty}(1-q_j)=+\infty\, .
$$
Together with \eqref{e:recursion}, this implies that $\lim\limits_{j\to \infty}p_{m_j}=0$.

It remains to prove \eqref{e:qj-diverge}.

Let $\Omega_j:=B_{m_{j+1}}\setminus \bigcup_{n\in I_{m_j}}\overline{B}(x_n,\tilde{r}_n)$, and define
$$
v(x):=\P_x\Big(X_{\tau_{\Omega_j}}\in \bigcup_{n\in I_{m_j}}\overline{B}(x_n,\tilde{r}_n)\Big)\, ,\quad x\in \R^d\, .
$$
Then $v$ is regular harmonic in $\Omega_j$. Note that
\begin{eqnarray*}
1-q_j&=&\inf_{x\in B_{m_{j+1}-2}\setminus B_{m_j}}\P_x(T_A <\tau_{B_{m_{j+1}}})\\
&\ge & \inf_{x\in B_{m_{j+1}-2}\setminus B_{m_j}} \P_x(X_{\tau_{\Omega_j}}\in \bigcup_{n\in I_{m_j}}\overline{B}(x_n,\tilde{r}_n))\\
&= & \inf_{x\in B_{m_{j+1}-2}\setminus B_{m_j}} v(x)\, .
\end{eqnarray*}
Hence, by \eqref{e:smaller-series}, it suffices to show that there exists $\wt{C}_1>0$ such that for all $j\in \N$ large enough, we have
\begin{equation}\label{e:inf-larger}
v(x)\ge \wt{C}_1 \sum_{n\in I_{m_j}}\frac{G(|x_n|)}{G(\tilde{r}_n)}\, \quad \textrm{for all }  x\in B_{m_{j+1}-2}\setminus B_{m_j}\, .
\end{equation}

Define
$$
u(x):=\sum_{n\in I_{m_j}}\P_x(T_{\overline{B}(x_n,\tilde{r}_n)}<\infty)\, .
$$
Then $u$ is regular harmonic in $\Omega_j$ and $u=1$ on each ball $\overline{B}(x_n,\tilde{r}_n)$, $n\in I_{m_j}$.

Let us show that there exists $\wt{C}_2=\wt{C}_2(\rho)>0$ such that
\begin{equation}\label{e:estimate-u-1}
u(x)\le \wt{C}_2\, ,\qquad \textrm{for all }x\in \bigcup_{n\in I_{m_j}} \overline{B}(x_n,\tilde{r}_n)\, .
\end{equation}

We are going to use a scaling-type argument. Let $a=\nu \rho^{-(m_j-1)}$ and let $X^a=(X^a_t)_{t\geq 0}$ be a L\'evy process with the L\'evy exponent $\psi^a$. Since $\P_{x}(T_{\overline{B}(x_n,\tilde{r}_n)}<\infty)=\P_{ax}(T^a_{\overline{B}(ax_n,a\tilde{r}_n)})$,  it is enough to show that
$$
	\sum\limits_{n\in I_{m_j}}\P_{y}(T^a_{\overline{B}(a x_n, a\tilde{r}_n)}<\infty)\leq \wt{C}_2,\ \qquad \text{ for all }\quad y\in \bigcup\limits_{n\in I_{m_j}}\overline{B}(a x_n, a\tilde{r}_n)
$$
for some constant $\wt{C}_2=\wt{C}_2(\rho)>0$.

For $n\in I_{m_j}$ set
\begin{align*}
    y_n & =a x_n & s_n& =a\tilde{r}_n\\  \tilde{s}_n&=G^a(s_n)^{-1/d} & \rho_n(y)&=\text{dist}(y,\overline{B}(y_n,s_n)) \, ,
\end{align*}
and note that $|y_n|=a|x_n|\le \nu \rho^{-(m_j-1)}\rho^{m_j}=\nu \rho$, $a|x_n|\ge \nu \ge 1$, $a r_n\le a|x_n| \le \nu\rho$, and $s_n=a\delta r_n \le \nu\rho$. In particular, $B(y_n,s_n)\subset B(0,2\nu\rho)$. Further, since $s_n=\delta(a r_n)$ and $a r_n\le \nu\rho$, we get from \eqref{eq:green_scale} that
$G^a(s_n)=G^a(\delta(a r_n))\le 2^{-d}(\delta(a r_n))^{-d}=2^{-d}s_n^{-d}$. Thus $\wt{s}_n=G^a(s_n)^{-1/d}\ge 2 s_n$.

Let $\mu_n$ be the equilibrium measure of $\overline{B}(y_n,s_n)$ with respect to the potential $G^a$. Then
$$
	\P_y(T^a_{\overline{B}(y_n,s_n)}<\infty)=(G^a\mu_n)(y)
	=\int\limits_{\overline{B}(y_n,s_n)}G^a(z-y)\mu_n(dz)\, ,\quad y\in \R^d\,.
$$	
Define the measure $\widetilde{\mu}_n$ by $\widetilde{\mu}_n(dy)=1_{\overline{B}(y_n,\tilde{s}_n)}(y)\,dy$\,.
We will prove that there is a constant $b>0$ such that for any $n\in I_{m_j}$
\begin{equation}\label{eq:todo2}
	G^a\widetilde{\mu}_n\geq b\, G^a\mu_n\quad \text{ on }\ \overline{B}(y_n,\widetilde{s}_n)^c\, .
\end{equation}	

Take $y\in\overline{B}(y_n,\widetilde{s}_n)^c$. Then
\begin{align*}
    G^a \widetilde{\mu}_n(y)&\geq  \int\limits_{B(y_n,\tilde{s}_n)\cap B(y,\rho_n)}G^a(z-y)\, dz\\
	&\geq G^a(\rho_n(y))\ |B(y_n,\tilde{s}_n)\cap B(y,\rho_n)|\\
	&\geq c_1\, G^a(\rho_n(y))\ |B(y_n,\tilde{s}_n)|\,,
\end{align*}
where in the last inequality we have used  \cite[Lemma 2.1]{AB} and the fact that $2 s_n \le \wt{s}_n$. On the other hand, by using Lemma \ref{l:capacity-estimate} in the last line,
\begin{align*}
    G^a \mu_n(y)&= \int\limits_{B(y_n,s_n)}G^a(z-y)\mu_n(dz)\\
	&\leq G^a(\rho_n(y))\mu_n(B(y_n,s_n))\\
	&=G^a(\rho_n(y))\cp^a(B(y_n,s_n))\\
	&\leq C_3 G^a(\rho_n(y)) G^a(s_n)^{-1}\,.
\end{align*}
Now we get (\ref{eq:todo2}) from the last two displays and the choice of $\tilde{s}_n$\,.

By (\ref{eq:lower_scaling_condition}), \eqref{eq:green_scale2}, (\ref{eq:cond_variant}) and (\ref{eq:choice_nu_delta})
\begin{align*}
	\frac{|y_j-y_k|^d}{(2\tilde{s}_k)^d}&\geq 2^{-d}a^d|x_j-x_k|^d G^a(s_k)=2^{-d}\psi^\star(a)|x_j-x_k|^d G(r_k)\\
    &\geq c_0 2^{-d}\psi^\star(a)\psi^\star(|x_k|^{-1})^{-1} =c_0 2^{-d}\frac{\psi^\star(a|x_k||x_k|^{-1})}{\psi^\star(|x_k|^{-1})}\\
	&	>c_0 C_L 2^{-d}\nu^\alpha=1,
\end{align*}
since $a|x_k|\geq \nu\geq 1$\,. Thus  $|y_j-y_k|> 2\max\{\tilde{s}_k,\tilde{s}_j\}$, which  implies  that the balls $\overline{B}(y_n,\tilde{s}_n)$ are disjoint.

Let $y\in \bigcup\limits_{n\in I_{m_j}}\overline{B}(y_n,s_n)$. Then there exists $n_0\in I$ so that $y\in \overline{B}(y_{n_0},s_{n_0})$ and $y\not\in \overline{B}(y_{n},\tilde{s}_{n})$ for all $n\in I_{m_j}\setminus\{n_0\}$\,.

Therefore, by (\ref{eq:todo2}) and Lemma \ref{lem:scaling} it follows that
\begin{align*}
	\sum\limits_{n\in I_{m_j}}\P_y(T_{\overline{B}(y_n,s_n)}<\infty)&\leq 1+\sum\limits_{\substack{n\in I_{m_j}\\ n\not=n_0}}(G^a\mu_n)(y)\leq 1+b^{-1}\sum\limits_{\substack{n\in I_{m_j}\\ n\not=n_0}}(G^a\widetilde{\mu}_n)(y)\\
    &\leq1+b^{-1}\int\limits_{B(0,2\nu\rho)}G^a(y)\,dy\\
    &\leq  1+ b^{-1}C_1\int\limits_{B(0,2\nu\rho)}\left(|y|^{2-d}\vee |y|^{\alpha-d}\right)\,dy=:\wt{C}_2\, .
\end{align*}
This proves \eqref{e:estimate-u-1}. Clearly, $v(x)= 1$ for $x\in \bigcup\limits_{n\in I_{m_j}} \overline{B}(x_n,\tilde{r}_n)$.

Now take $x\in B_{m_{j+1}}^c$. Then
$$
|x-x_n|\ge  |x|-|x_n|\ge (\rho^{m_{j+1}}-\rho^{m_j})=(\rho^{m_0}-1)\rho^{m_j}\ge (\rho^{m_0}-1) |x_n|\, ,
$$
implying that
\begin{eqnarray*}
u(x)&\le &\sum_{n\in I_{m_j}}\frac{G(|x-x_n|)}{G(\tilde{r}_n)}\le \sum_{n\in I_{m_j}}\frac{G\big((\rho^{m_0}-1)|x_n|\big)}{G(\tilde{r}_n)}\\
&\leq& 4C_G^2\left(\frac{\rho^{m_0-1}+\rho}{\rho^{m_0}-1}\right)^{d-2}\sum_{n\in I_{m_j}}\frac{G((\rho^{m_0-1}+\rho)|x_n|)}{G(\tilde{r}_n)}\\
&\leq & 4C_G^2\left(\frac{2\rho}{\rho^2-1}\right)^{d-2}\sum_{n\in I_{m_j}}\frac{G((r^{m_0-1}+\rho)|x_n|)}{G(\tilde{r}_n)}\\
&\leq &\frac{C_E}{2}\sum_{n\in I_{m_j}}\frac{G((r^{m_0-1}+\rho)|x_n|)}{G(\tilde{r}_n)}\, ,
\end{eqnarray*}
where the third inequality we have used Lemma \ref{lem:green_est}, in the fourth that $m_0\geq 2$ and in the last one (\ref{e:choice-rho}). Clearly, $v(x)=0$ for $x\in B_{m_{j+1}}^c$. By the maximum principle (Lemma \ref{lem:maxpr}) it follows that
\begin{equation}\label{e:max-v-u}
\wt{C}_2 v(x)\ge u(x)- \frac{C_E}{2}  \sum_{n\in I_{m_j}}\frac{G\big((\rho^{m_0-1}+\rho)|x_n|\big)}{G(\tilde{r}_n)}\, ,\quad \textrm{for all }x\in \R^d\, .
\end{equation}

Take $x\in B_{m_{j+1}-2}\setminus B_{m_j}$. Then for $x_n$ with index in $I_{m_j}$,
$$
|x-x_n|\le |x|+|x_n|\le \rho^{m_{j+1}-2}+\rho^{m_j}=(\rho^{m_0-1}+\rho)\rho^{m_j-1}\le (\rho^{m_0-1}+\rho)|x_n|\, .
$$
By the left-hand side inequality of Lemma \ref{l:equilibrium-potential} we have
\begin{equation}\label{e:estimate-u-2}
u(x)\ge C_E \sum_{n\in I_{m_j}}\frac{G(|x-x_n|)}{G(\tilde{r}_n)} \ge C_E \sum_{n\in I_{m_j}} \frac{G\big((\rho^{m_0-1}+\rho)|x_n|\big)}{G(\tilde{r}_n)}\, .
\end{equation}
Combining \eqref{e:max-v-u} and \eqref{e:estimate-u-2} gives
$$
v(x)\ge \frac{C_E}{2\wt{C}_2} \sum_{n\in I_{m_j}}\frac{G\big((\rho^{m_0-1}+\rho)|x_n|\big)}{G(\tilde{r}_n)}\, ,\quad x\in B_{m_{j+1}-2}\setminus B_{m_j}\, .
$$
Finally, by  (\ref{eq:lower_scaling_condition}) and Lemma \ref{lem:green_est}, $G\big((\rho^{m_0-1}+\rho)|x_n|\big)\ge G\big(\rho^{m_0}|x_n|\big)\ge \frac{\rho^{m_0(2-d)}}{4C_G^2}\, G(|x_n|)$. Hence,
$$
v(x)\ge \wt{C}_1 \sum_{n\in I_{m_j}}\frac{G(|c_n|)}{G(\tilde{r}_n)}\, ,\quad x\in B_{m_{j+1}-2}\setminus B_{m_j}\, ,
$$
with $\wt{C}_1:=\frac{C_E\rho^{m_0(2-d)}}{8\wt{C}_2 C_G^2}$. \qed

\section{Proof of Theorem \ref{t:theorem-2}}\label{s:proof-2}

\begin{Lem}\label{l:series-integral}
If the family of balls  $\{\overline{B}(x_n,r_n)\}_{n\ge 1}$ is regularly located, then
$$
\sum_{n=1}^{\infty} \frac{G(|x_n|)}{G(\phi(|x_n|))}<\infty \quad \textrm{if and only if } \quad \int\limits_1^\infty \frac{r^{d-1}G(r)}{G(\phi(r))}\, dr <\infty\, .
$$
\end{Lem}
\proof
Recall that $R>0$ is the constant from the uniform density condition.
 Set $A_k=\{n\in \N\colon 2Rk\leq |x_n|<2R(k+1)\}$ for $k\in \N$\,. By comparing measures of the sets and using the $\epsilon$ separation condition it follows that the number of balls with centers in $A_k$ is at most $\frac{(2R(k+1)+\epsilon)^d-(2Rk-\epsilon)^d)}{\epsilon^d}$ for $k$ large enough so that $r_n<\epsilon$. Hence, for some constant $c_1=c_1(\epsilon,d,R)>0$, there are at most $c_1(2kR)^{d-1}$ balls with center in $A_k$ for all $k\in \N$\,.
	
In a similar fashion, by using uniform density, we deduce that, for some constant $c_2=c_2(\epsilon,d,R)>0$,  there are at least $c_2(2kR)^{d-1}$ balls with center in $A_k$\,.
	
Since both $G$ and $\phi$ are decreasing we obtain
	\begin{align*}
		\sum\limits_{n=1}^\infty \frac{G(|x_n|)}{G(r_n)}&=\sum\limits_{n=1}^\infty \frac{G(|x_n|)}{G(\phi(|x_n|))} \geq \sum\limits_{k=1}^\infty \frac{G(2R(k+1))}{G(\phi(2R(k+1)))}c_2(2kR)^{d-1}\\
		&\geq c_3 \sum\limits_{k=1}^\infty \int\limits_{2R(k+1)}^{2R(k+2)}\frac{G(r)}{G(\phi(r))}r^{d-1}\,dr\geq c_4\int\limits_1^\infty \frac{r^{d-1}G(r)}{G(\phi(r))}\, dr\, .
	\end{align*}	
The other bound can be proved similarly.
\qed

By Theorem \ref{t:theorem-1} if \eqref{e:thm-2} holds then the collection of balls is avoidable. In the other direction, we effectively show that if the regularly spaced collection of balls is avoidable then the separation condition \eqref{e:thm-1-2} is automatically satisfied.
We will need to count the balls in annuli of the form $A(x,r,2 r):=\{y\in \R^d:\, r\le |y-x|\le 2 r\}$, $r>0$.
From the separation condition it follows that there exists $N_1=N_1(d)\in \N$ such that for all $x\in \R^d$ and all $r>0$ the number of balls $\overline{B}(x_n,r_n)$ with centers in $A(x,r,2 r)$ is bounded from above by $N_1 r^d$. Similarly, the uniform density condition implies that there exists $N_2=N_2(d)\in \N$ such that for all $r>0$ large enough the number of balls $\overline{B}(x_n,r_n)$ with centers in $A(0,r,2 r)$ is bounded from below by $N_2 r^d$.

Recall that $A=\bigcup_{n\ge 1}\overline{B}(x_n,r_n)$.

\begin{Lem}\label{l:reg-loc}
Suppose that the family of balls  $\{\overline{B}(x_n,r_n)\}_{n\ge 1}$ is regularly located. If the function $r\mapsto r^d G(r)/G(\phi(r))$ is unbounded (at infinity), then $A$ is unavoidable.
\end{Lem}
\proof Assume that $A$ is avoidable and let $\theta:=\P_0(T_A=\infty)>0$. Fix
$$
a:=\left(\frac{2C_E}{\theta}\right)^{\alpha}\vee 4\, .
$$
Since $r^d G(r)/G(\phi(r))$ is unbounded, there exists a sequence $\{R_j\}_{j\ge 1}$ such that $(R_j)^{d}G(R_j)/G(\phi(R_j))\to \infty$ as $j\to \infty$, and additionally satisfies  $R_{j+1}\ge 4a R_j$ and
\begin{equation}\label{e:aux2}
\frac{G\left(\frac{R_{j+1}}{2}\right)}{G(2aR_j)} \le \frac{C_E}{2}\frac{N_2}{N_1}\, .
\end{equation}
For $j\ge 1$ let $B_j:=B(0,R_j)$,
\begin{eqnarray*}
p_j&:=&\P_0(\tau_{B_j}<T_A)\, ,\\
q_j&:=&\sup_{x\in B(0,aR_j)\setminus B_j}\P_x(\tau_{B_{j+1}}<T_A)\, ,\\
t_j&:=&\P_0\Big(X_{\tau_{B_j}}\in B(0,aR_j)\, \Big|\,  \tau_{B_j}<T_A\Big)\, .
\end{eqnarray*}
By Lemma \ref{lem:recurrence}, for all $j\ge 1$,
$$
p_{j+1}\le p_j(q_j t_j+1-t_j)\, .
$$
Moreover, by Lemma \ref{l:long-jump},
$$
t_j\ge \frac{p_j-\P_0(X_{\tau_{B_j}}\notin B(0,aR_j))}{p_j}\ge \frac{\theta-\frac{\theta}{2}}{p_j}\ge \frac{\theta}{2}\, .
$$
Let $I_j:=\{n\ge 1:\, R_j\le |x_n|\le 2R_j\}$ be the set of indices of balls with centers in the annulus $A(0,R_j, 2R_j)$. Define $A_j:=\bigcup_{n\in I_j}\overline{B}(x_n,r_n)$. We claim that there exists $\delta >0$ such that
\begin{equation}\label{e:lower-bound}
\P_x(T_{A_j}<\tau_{B_{j+1}})\ge \delta\, ,\qquad \textrm{for all } x\in B(0,aR_j)\setminus B_j\, .
\end{equation}
This will imply that for all $j\ge 1$,
$$
1-q_j=\inf_{x\in B(0,aR_j)\setminus B_j}\P_x(T_A<\tau_{B_{j+1}})\ge \inf_{x\in B(0,aR_j)\setminus B_j}\P_x(T_{A_j}<\tau_{B_{j+1}})\ge \delta\, .
$$
In the same way as in the proof of Theorem \ref{t:theorem-1}(b) we conclude that $\lim_{j\to \infty}p_j=0$, thus deriving a contradiction. Hence, it remains to show \eqref{e:lower-bound}. We first modify the function $\phi$ in the following way. Let $\wt{\phi}(r)=\phi(2R_j)$ if $r\in [R_j,2R_j]$ for some $j\ge 1$, and $\wt{\phi}(r)=\phi(r)$ otherwise. Take new balls $\overline{B}(x_n,\wt{\phi}(|x_n|))$ and let $\wt{A}_j:=\bigcup_{n\in I_j}\overline{B}(x_n,\wt{\phi}(|x_n|))$. Since $\wt{A}_j\subset A_j$, proving \eqref{e:lower-bound} for $\wt{A}_j$ will suffice. We keep using notation $\phi$ for the new $\wt{\phi}$. Then for $n\in I_j$ we have that $r_n=\phi(|x_n|)=\phi(2R_j)=\phi(R_j)=:\phi_j$.

Define $u:\R^d \to [0,\infty)$ by
$$
u(x):=\sum_{n\in I_j} \P_x(T_{\overline{B}(x_n,r_n)}<\infty)\, .
$$
Then $u$ is harmonic in $\Omega_j:=B_{j+1}\setminus A_j$.

Let $n_0\in I_j$ and $x\in B(x_{n_0}, r_{n_0})$. By the separation condition, for any $k\ge 1$, there are at most $N_1 2^{kd} $ balls with centers $x_n$ such that $\epsilon 2^{k-1}\le |x-x_n|\le \epsilon 2^k$. For such $x_n$ we have by Lemma \ref{l:equilibrium-potential}
$$
\P_x(T_{\overline{B}(x_n,r_n)}<\infty)\le \frac{G(|x-x_n|)}{G(\phi(|x_n|))}\le \frac{G(\epsilon 2^{k-1})}{G(\phi_j)} \, .
$$
The at most $N_1 2^{kd}$ such balls contribute in the sum (for $u$)  at most $N_1 2^{kd} G(\epsilon 2^{k-1})/G(\phi_j)$. It suffices to count only such balls for which $\epsilon 2^k \le 6 R_j$ (there are none for which this is not true). Now we estimate
$$
\sum_{k:\, \epsilon 2^k\le 6R_j} 2^{kd} G(\epsilon 2^{k-1})\, .
$$
Let $k_0\in \N$ be such that $\epsilon 2^{k_0-1}\le 6R_j <\epsilon 2^{k_0}$. Then by using \eqref{e:green-est} in the second and the sixth inequality, \eqref{eq:lower_scaling_condition} in the third, and \eqref{eq:square} in the fifth, we get
\begin{eqnarray*}
\sum_{k:\, \epsilon 2^k\le 6R_j} 2^{kd} G(\epsilon 2^{k-1})&\le & \sum_{k=1}^{k_0} 2^{kd} G(\epsilon 2^{k-1})
\le\frac{2^d C_G}{\epsilon^d}\sum_{k=1}^{k_0}  \frac{1}{\psi^{\star}((\epsilon 2^{k-1})^{-1})}\\
&\le & \frac{2^d C_G}{\epsilon^d C_L}\frac{1}{\psi^{\star}((\epsilon 2^{k_0-1})^{-1})}\sum_{k=1}^{k_0}2^{(k-k_0)\alpha}\\
&\le & \frac{2^d C_G}{\epsilon^d C_L(1-2^{-\alpha})}\frac{1}{\psi^{\star}((6R_j)^{-1})}\\
&\le & \frac{2(1+6^2)2^d C_G }{\epsilon^d C_L(1-2^{-\alpha})}\frac{1}{\psi^{\star}((R_j)^{-1})}\\
&\le & \frac{74\, 2^d C_G^2 }{\epsilon^d C_L(1-2^{-\alpha})}\, G(R_j) R_j^d =c_1 G(R_j) R_j^d
\end{eqnarray*}
with $c_1=c_1(d,\psi)>0$.
Hence,
$$
u(x) \le  1+c_1 \frac{G(R_j) R_j^d}{G(\phi(R_j))}\le c_2 \frac{G(R_j) R_j^d}{G(\phi(R_j))}\, ,
$$
by the assumption that $r^dG(r)/G(\phi(r))$ is unbounded. Thus we have proved that
\begin{equation}\label{e:u-estimate1b}
u(x) \le c_2\frac{G(R_j) R_j^d}{G(\phi_j)}\, ,\qquad x\in A_j\, .
\end{equation}

Let $x\in B(0,aR_j)\setminus B_j$. Then $|x-x_n|\le |x|+|x_n|\le aR_j+R_j\le 2a R_j$. Thus
\begin{equation}\label{e:u-estimate2b}
u(x) \ge  C_E \sum_{n\in I_j}\frac{G(|x-x_n|)}{G(\phi_j)} \ge C_E \sum_{n\in I_j}\frac{G(2 a R_j)}{G(\phi_j)} \ge C_E \frac{G(2 a R_j)}{G(\phi_j)} N_2 R_j^d\, .
\end{equation}

Let $x\in B_{j+1}^c$. Then $|x-x_n|\ge |x|-|x_n|\ge R_{j+1}-2R_j\ge R_{j+1}/2$. Hence,
\begin{eqnarray}
u(x) &\le &  \sum_{n\in I_j}\frac{G(|x-x_n|)}{G(\phi_j)} \le  \sum_{n\in I_j}\frac{G(R_{j+1}/2)}{G(\phi_j)}\le \frac{G(R_{j+1}/2)}{G(\phi_j)} N_1 R_j^d \nonumber \\
& \le &\frac{N_1 R_j^d}{G(\phi_j)}\, \frac{C_E}{2}\frac{N_2}{N_1} G(2a R_j) = \frac{C_E}{2} \frac{G(2a R_j)}{G(\phi_j)} N_2 R_j^d\, , \label{e:u-estimate3b}
\end{eqnarray}
where in the second line we used \eqref{e:aux2}.

The maximum principle together with \eqref{e:u-estimate1b} and \eqref{e:u-estimate3b} implies that
\begin{equation}\label{e:u-estimate4}
c_2 \frac{G(R_j)R_j^d}{G(\phi_j)}\,  \P_x(T_{A_j}<\tau_{B_{j+1}})\ge u(x)-\frac{N_2 C_E}{2} \frac{G(2a R_j)}{G(\phi_j)} R_j^d\, ,\qquad x\in \Omega_j\, .
\end{equation}
Note that by \eqref{e:green-est} and \eqref{eq:square},
$$
G(2a R_j)\ge \frac{1}{2^{d+1}  (1+4a^2)C_G^2}\, G(R_j)\, .
$$
By using \eqref{e:u-estimate2b} for $x\in B(0,aR_j)\setminus B_j$ and the above estimate, we get
$$
c_2 \frac{G(R_j)R_j^d}{G(\phi_j)}  \P_x(T_{A_j}<\tau_{B_{j+1}})\ge \frac{N_2 C_E}{2} \frac{G(2a R_j)}{G(\phi_j)} R_j^d
\ge \frac{N_2 C_E}{2^{d+2}  (1+4a^2)C_G^2} \frac{G( R_j)R_j^d}{G(\phi_j)} \, .
$$
Hence
$$
\P_x(T_{A_j}<\tau_{B_{j+1}})\ge \frac{N_2 C_E}{2^{d+2}  (1+4a^2)C_G^2 c_2}\, .
$$
This proves \eqref{e:lower-bound}. \qed

\emph{Proof of Theorem \ref{t:theorem-2}.}
Assume that the integral in \eqref{e:thm-2} converges. By Lemma \ref{l:series-integral},  $\sum_{n\ge 1} \frac{G(|x_n|)}{G(\phi(|x_n|))}<\infty$. Proposition \ref{p:sufficient-cond} now implies that the collection of balls is avoidable.

To show the converse, assume that the (regularly located) collection of balls is avoidable. It follows from Lemma \ref{l:reg-loc} that $r\mapsto r^d G(r)/G(\phi(r))$ is bounded (at infinity) implying that
$$
\inf_{n\ge 1} \frac{G(\phi(|x_n|))}{|x_n|^d G(|x_n|)}>0\, .
$$
By \eqref{eq:comp} and \eqref{e:green-est},
$$
\psi(|x_n|^{-1})\asymp \psi^{\star}(|x_n|^{-1})\asymp |x_n|^{-d}G(|x_n|)^{-1}\, .
$$
Together with $\epsilon$ separation this gives
$$
\inf\limits_{m\not= n} |x_m-x_n|^d\, \psi(|x_n|^{-1})\, G(\phi(|x_n|))\ge (2\epsilon)^d \inf_{n\ge 1}\frac{G(\phi(|x_n|))}{|x_n|^d G(|x_n|)} >0\, .
$$
Thus, the separation condition \eqref{e:thm-1-2} is satisfied and the claim follows from Theorem \ref{t:theorem-1}(b). \qed

\section{Poissonian collections of balls}\label{s:poisson}

In this section we consider collections of randomly located closed balls with centers coming from a Poisson point process. Since we closely follow the arguments from \cite{CDOC}, most of the proofs are omitted. We precisely spell out the assumptions and give proofs when some (non-trivial) changes are needed.

We consider a Poisson point process on $\R^d$ with mean measure $\mu(dx)=\mu(x)\, dx$ absolutely continuous with respect to the Lebesgue measure and an independent isotropic unimodal L\'evy process $X$ with the characteristic exponent $\psi$ satisfying the weak lower scaling condition \eqref{eq:lower_scaling_condition}. In this section we will assume that $X$ is a subordinate Brownian motion. Let $G$ denote the Green function of $X$. We impose the following assumptions on the radius function $\phi$ and the density $\mu$: There exists a constant $C_P>1$ such that
    \begin{equation}\label{e:poisson-condition-1}
    \begin{array}{ll}
    \left\{\begin{array}{l}C_P^{-1}\phi(x)\le \phi(y)\le C_P \phi(x)\\
    C_P^{-1}\mu(x)\le \mu(y)\le C_P \mu(x)\end{array}\right. & \textrm{for } y\in B\left(x,\frac{|x|}{2}\right)\, ,
    \end{array}
    \end{equation}
    \begin{equation}\label{e:poisson-condition-2}
    \phi(x)\le \frac12 |x|\, ,\quad \textrm{for all }x\in \R^d\, ,
    \end{equation}
    \begin{equation}\label{e:poisson-condition-3}
       \psi^\star(|x|^{-1})^{-1} G(\phi(x))^{-1}\mu(x)\le C_P\, ,\quad |x|\to \infty\, .
    \end{equation}

Let $\sP$ be a realization of points from the Poisson point process and let
\begin{equation}\label{e:random-balls}
A_{\sP}=\bigcup_{x\in \sP} \overline{B}(x,\phi(x))\, .
\end{equation}
The collection of closed balls $A_{\sP}$ is said to be avoidable if there exists a point $x\in \R^d$ such that $\P_x(T_{A_{\sP}}<\infty)<1$. Note that  $A_{\sP}$ is avoidable if and only if it is minimally thin at infinity. The latter condition apparently does not depend on the starting point $x$.

Following \cite{CDOC} we say that we have \emph{percolation L\'evy process} if there is a positive probability that the realization of point from the Poisson point process results in an avoidable collection of balls.

\begin{Thm}\label{t:theorem-poisson}
Suppose that \eqref{e:poisson-condition-1}--\eqref{e:poisson-condition-3} hold and that $X$ is a subordinate Brownian motion with the characteristic exponent $\psi$ satisfying the weak lower scaling condition \eqref{eq:lower_scaling_condition}. Percolation L\'evy process occurs if and only if
\begin{equation}\label{e:thm-poisson}
\int_{|x|>1}\frac{G(x)}{G(\phi(x))}\, \mu(x)\, dx <\infty\, .
\end{equation}
Moreover, in case percolation L\'evy process occurs, the random collection of balls $A_{\sP}$ is avoidable with probability one.
\end{Thm}

In order to prove this theorem we use the characterization of minimal thinness at infinity given in Theorem \ref{t:aikawa-criterion}. Let $\{Q_m\}_{m\ge 1}$ be the Whitney decomposition described in the introduction and let
$$
W(A_{\sP}, \infty)=\sum_{m\ge 1}G(\mathrm{diam}(Q_m))\, \mathrm{Cap}(A_{\sP}\cap Q_m)\, .
$$
Then $A_{\sP}$ is avoidable if and only if $W(A_{\sP}, \infty)<\infty$.

The following lemma is proved exactly in the same way as \cite[Lemma 2, Lemma 4]{CDOC}.
\begin{Lem}\label{l:poisson-0-1}
\begin{itemize}
    \item[(i)] $\P(A_{\sP} \textrm{ is minimally thin at infinity} ) \in \{0,1\}$.
    \item[(ii)] $\E[W(A_{\sP}, \infty)]<\infty$ if and only if $W(A_{\sP}, \infty)<\infty $ $\P$-a.s.
\end{itemize}
\end{Lem}
The next result is the key part of the proof of the theorem.

\begin{Lem}\label{l:poisson-key}
There exists a constant $C_4>1$ such that for any Whitney cube $Q$ and any point $x\in Q$
    \begin{equation}\label{e:poisson-key}
    C_4^{-1}  \frac{\mu(Q)}{G(\phi(x))}\le \E[\cp(A_{\sP}\cap Q)]\le C_4  \frac{\mu(Q)}{G(\phi(x))}\, .
    \end{equation}
\end{Lem}
We defer the proof of the lemma until the end of the section and outline the proof of theorem based on Lemmas \ref{l:poisson-0-1} and \ref{l:poisson-key}.

\emph{Proof of Theorem \ref{t:theorem-poisson}.} Note that
$$
\E[W(A_{\sP}, \infty)]=\sum_{m=1}^{\infty}G(\mathrm{diam}(Q_m)) \, \E[\cp(A_{\sP}\cap Q_m)]\, .
$$
By Lemma \ref{l:poisson-key}, for any point $x_m\in Q_m$,  $\E[\cp(A_{\sP}\cap Q_m)]\asymp \frac{\mu(Q_m)}{G(\phi(x_m))}$. Since $\mathrm{diam}(Q_m)\asymp |x_m|$, we also have that $G(\mathrm{diam}(Q_m))\asymp G(|x_m|)$. Hence the sum above is comparable with
$$
\sum_{m=1}^{\infty} \frac{G(|x_m|)}{G(\phi(x_m))}\, \mu(Q_m)\, .
$$
It follows from \eqref{e:poisson-condition-1} that both $\phi$ and $\mu$ are approximately constant on each Whitney cube. This implies that the last sum is convergent if and only if the integral in \eqref{e:thm-poisson} is convergent. Hence, $\E[W(A_{\sP}, \infty)]<\infty$ if and only if the integral in \eqref{e:thm-poisson} converges. The first claim of the theorem now follows form Lemma \ref{l:poisson-0-1}(ii), while the second one from Lemma \ref{l:poisson-0-1}(i). \qed

\emph{Proof of Lemma \ref{l:poisson-key}.} The right-hand side inequality is proved exactly in the same way as in the proof of \cite[Lemma 3]{CDOC}. To prove the left-hand side inequality we follow \cite{CDOC} and use the super-additivity property of capacity due to Aikawa and Borichev, \cite[Theorem 3]{AB}, which in our case reads as follows: For $r>0$ and $a>0$ the radius $\eta^a(r)$ is chosen so that $\cp^a(B(0,r))=|B(0,\eta^a(r))|=\sigma_d \eta^a(r)^d$ ($\sigma_d$ is the volume of the unit ball). By Lemma \ref{l:capacity-estimate},
$$
C_3^{-1}G^a(r)^{-1}\le \cp^a (B(0,r)) \le C_3 G^a(r)^{-1}\, .
$$
This gives the following estimate for $\eta^a(r)$:
$$
\wt{C}_3^{-1}G^a(r)^{-1/d}\le \eta^a(r) \le \wt{C}_3G^a(r)^{-1/d},
$$
for $\wt{C}_3=\wt{C}_3(d, C_3)=\wt{C}_3(d,C_L, \alpha)\geq 1$.
Let $r_0>0$. Choose $a=a(r_0)>0$ small enough so that
\begin{equation}\label{eq:aikawa-borichew1}
\wt{C}_3^{-1}G^a(ar)^{-1/d}\ge 2ar\quad \text{ for all }\ \ 0<r\leq r_0\,.
\end{equation}

To see that this is possible, first note that
$$
G^a(a r)^{1/d}=a^{-1}\psi^{\star}(a)^{1/d} G(r)^{1/d}\ .
$$
Then it follows for $a\leq r_0^{-1}\le r^{-1}$ that
\begin{align*}
\wt{C}_3^{-1}G^a(ar)^{-1/d}&= \wt{C}_3^{-1}  a\psi^\star(a)^{-1/d}G(r)^{-1/d}\geq C_G^{-1/d}\wt{C}_3^{-1}ar\left(\frac{\psi^\star(r^{-1})}{\psi^\star(a)}\right)^{1/d}\\
&\geq C_L^{2/d}C_G^{-1/d}\wt{C}_3^{-1}ar (ar)^{-\alpha}
\end{align*}
 and so it is enough to choose $a>0$ small enough so that
\begin{equation}\label{eq:a_cond}
	C_L^{2/d}C_G^{-1/d}\wt{C}_3^{-1}(ar_0)^{-\alpha}\geq 2
\end{equation}
i.e.
\begin{equation}\label{eq:choice_of_a}
	a\leq \left(\frac{2C_G^{1/d}\wt{C}_3}{C_L^{2/d}}\right)^{-1/\alpha}r_0^{-1}\,.
\end{equation}
Note that $a\le r_0^{-1}$.
Hence, $(\eta^a)^*(ar):=\max(\eta^a(ar),2ar)=\eta^a(ar)$ for $0<r\leq r_0$.

Let $F:=\bigcup B(y_k,a\rho_k)\subset B(x,1)$ for some $x\in \R^d$, where $\rho_k\le r_0$ and the larger balls
$B(y_k, \wt{C}_3G^a(a\rho_k)^{-1/d})$ are disjoint (note that the latter balls are larger because of (\ref{eq:aikawa-borichew1})). Then $B(x_k, (\eta^a)^*(a\rho_k))$ are disjoint, hence by \cite[Theorem 3]{AB}
\begin{equation}\label{e:aikawa-borichev}
\cp^a(F)\ge c_1\sum_k \cp^a(B(y_k,a\rho_k))\, ,
\end{equation}
for some constant $c_1>0$.

Let $\phi_0=\min_{x\in Q}\phi(x)$. By \eqref{e:poisson-condition-1}, $\phi_0\asymp \phi(x)$ for all $x\in Q$. It suffices to consider only balls with centers in $Q$ and assume that all such balls are of radius $\phi_0$ (this decreases the capacity of $\sP \cap Q$).
Choose (cf. (\ref{eq:choice_of_a}) with $r_0=\phi_0$)
$$
a=\min\left\{\frac{1}{\ell(Q)\sqrt{d}},\left(\frac{2C_G^{1/d}\wt{C}_3}{C_L^{2/d}}\right)^{-1/\alpha}\phi_0^{-1}\right\}\,,
$$
where $\ell(Q)$ is the sidelength  of $Q$.

Note that $a\leq \phi_0^{-1}$ and that \eqref{e:poisson-condition-2} implies $a\ge c_2/\ell(Q)$.

Let
$$
N=\lfloor (8C_G^{1/d} \wt{C}_3)^{-1} \ell(Q) \, a G^a(a\phi_0)^{1/d} \rfloor \, .
$$

The cube $Q$ is divided into $N^d$ smaller sub-cubes with sidelength $\ell(Q)/N$, and a typical sub-cube is denoted by $Q'$. Let $Q''\subset Q'$ be a cube concentric to $Q'$ with sidelength $\ell(Q)/(4N)$. If $Q''$ contains points from the realization of $\sP$, choose one such point. This gives points $x_1,x_2, \dots, x_M$, where $M\le N^d$ is random. By the choice of $a$ and $N$, $B(x_k,\phi_0)\subset Q'$ where $Q'$ is the sub-cube that contains $x_k$. Indeed, since
$$
aG^a(a\phi_0)^{1/d}=\psi^{\star}(a)^{1/d} G(\phi_0)^{1/d}
$$
we get
\begin{eqnarray*}
N&\le & (8C_G^{1/d} \wt{C}_3)^{-1} \ell(Q)\psi^{\star}(a)^{1/d}G(\phi_0)^{1/d} \\
&\le & (8C_G^{1/d} \wt{C}_3)^{-1} \ell(Q)\psi^{\star}(a)^{1/d} C_G^{1/d}\phi_0^{-1}\psi^{\star}(\phi_0^{-1})^{-1/d}\\
&\le & 8^{-1}\ell(Q)\phi_0^{-1} \left(\frac{\psi^{\star}(a)}{\psi^{\star}(\phi_0^{-1})}\right)^{1/d}\le
 8^{-1}\ell(Q)\phi_0^{-1}\, ,
\end{eqnarray*}
where the last inequality follows from $a\le \phi_0^{-1}$.
This implies that  $\phi_0\le \frac{\ell(Q)}{8N}$.
Set
$$
A_{\sP, Q}=\bigcup_{k=1}^M \overline{B}(x_k, \phi_0)\, .
$$
Then $\cp(A_{\sP}\cap Q)\ge \cp(A_{\sP, Q})$.
To estimate $\cp(A_{\sP, Q})$ we scale the cube $Q$ by factor $a$. Since $\ell(aQ)\le 1/\sqrt{d}$, the scaled cube $aQ$ is inside a ball of radius 1.
By choice of $N$ we have that
$$
 8 C_G^{1/d} \wt{C}_3G^a(a\phi_0)^{-1/d}\le \frac{a\ell(Q)}{N}\, ,
$$
which implies that balls $B(a x_k, \wt{C}_3 G^a(a\phi_0)^{-1/d})$ are disjoint. The scaling relation for capacity (Lemma \ref{l:capacity-scaling}), \eqref{e:aikawa-borichev}, and Lemma \ref{l:capacity-estimate} imply that
\begin{eqnarray*}
\cp(A_{\sP, Q})&=&\psi^{\star}(a)a^{-d}\cp^a(a A_{\sP, Q})\\
&\ge & c_1 \psi^{\star}(a)a^{-d} M \cp^a(B(0, a\phi_0))\\
&=& c_1 M \cp(B(0,\phi_0))\ge \frac{c_3}{G(\phi_0)} \, M\, .
\end{eqnarray*}
Hence,
\begin{equation}\label{e:poisson-aux}
\E[\cp(A_{\sP}\cap Q)]\ge \frac{c_3}{G(\phi_0)}\, \E M\, .
\end{equation}
The probability that $Q''$ contains a Poisson point is equal to
$$
1-\P(\sP \cap Q''=\emptyset )=1-e^{-\mu(Q'')}\, .
$$
To estimate $\mu(Q'')$, let $x\in Q''$ be the center of $Q''$. Then
$$
\mu(Q'')\asymp \mu(x)|Q''| \asymp \mu(x)\left(\frac{\ell(Q)}{N}\right)^d\, .
$$

By choice of $N$ and using $a\ge c_2/\ell(Q)$ we get
\begin{eqnarray*}
\left(\frac{\ell(Q)}{N}\right)^d &\le  &c_4 a^{-d}G^a(a\phi_0)^{-1} =c_4\psi^\star(a)^{-1}G(\phi_0)^{-1}\\
&\leq & c_5\psi^\star(c_2/\ell(Q))^{-1}G(\phi_0)^{-1}\,.
\end{eqnarray*}
Thus, by (\ref{eq:square-left}) we get
$$
\mu(Q'')\le c_6 \psi^\star(\ell(Q)^{-1})^{-1}G(\phi_0)^{-1}\mu(x)\le c_{7} \psi^{\star}(|x|^{-1})^{-1} G(\phi(x))^{-1}\mu(x)\le c_{7} C_P
$$
by \eqref{e:poisson-condition-3}. Finally, since $1-e^{-x}\ge c_{8}x$ for $x\in [0, c_{7} C_P]$,
$$
\E M=\sum_{Q''\subset Q}\left(1-e^{-\mu(Q'')}\right)\ge c_{8}\sum_{Q''\subset Q}\mu(Q'')\ge c_{9}\mu(Q)\, .
$$
Together with \eqref{e:poisson-aux} this gives that $\E[\cp(A_{\sP}\cap Q)]\ge c_{9}\frac{\mu(Q)}{G(\phi(x))}$. \qed

\section{Case $d\le 2$}\label{s:d-le-2}
Let $f:(0,\infty)\to (0,\infty)$ be a Bernstein function with the representation
$$
f(\lambda)=\int_{[0,\infty)}(1-e^{-\lambda t})\, m(dt)\, ,\quad \lambda>0\, .
$$
If the L\'evy measure $m(dt)$ has a completely monotone density $m(t)$, then $f$ is called a complete Bernstein function. The function $f$ serves as the Laplace exponent of a subordinator (i.e.~a nonnegative L\'evy process) $S=(S_t)_{t\ge 0}$ -- the law of $S_t$ is characterized by $\E[e^{-\lambda S_t}]=e^{-t f(\lambda)}$. Let $W=(W_t,\P_x)$ be an independent $d$-dimensional Brownian motion. The process $X=(X_t)_{t\ge 0}$ defined by $X_t=W(S_t)$ is called a subordinate Brownian motion. It is an isotropic unimodal L\'evy process with the characteristic exponent $\psi$ given by $\psi(x)=f(|x|^2)$. Note that $\psi(|x|)$ is already increasing and hence equal to $\psi^{\star}$.

Following \cite{KSV8} we will assume the next two hypotheses:

\noindent
{\bf (H1):}
There exist constants $0<\delta_1\le \delta_2 <1$ and $a_1, a_2>0$  such that
\begin{equation}\label{e:H1}
a_1\lambda^{\delta_1} f(t) \le f(\lambda t) \le a_2 \lambda^{\delta_2} f(t), \quad \lambda \ge 1, t \ge 1\, .
\end{equation}
\noindent
{\bf (H2):}
There exist constants $0<\delta_3\le \delta_4 <1$ and $a_3, a_4>0$  such that
\begin{equation}\label{e:H2}
a_3\lambda^{\delta_4} f(t) \le f(\lambda t) \le a_4 \lambda^{\delta_3} f(t),
\quad \lambda \le 1, t \le 1\, .
\end{equation}
By \cite[(2.11)]{KSV8} we have that under {\bf (H1)} and {\bf (H2)} there are constants $a_5$ and $a_6$ such that
$$
a_5 \lambda^{\delta_1 \wedge \delta_3}\le \frac{f(\lambda r)}{f(r)} \le a_6 \lambda^{\delta_2\vee \delta_4}\, ,\quad \lambda\ge 1, r>0\, .
$$
In terms of the characteristic exponent $\psi$ the above reads
$$
a_5 \lambda^{2(\delta_1\wedge \delta_3)} \psi(r)\le \psi(\lambda r) \le a_6 \lambda^{2(\delta_2\vee \delta_4)} \psi(r)\, ,\quad \lambda\ge 1, r>0\, .
$$
If we define $C_L:=a_5$ and $\alpha:=2(\delta_1 \wedge \delta_3)$, we see that \eqref{eq:lower_scaling_condition} holds true. Note that $\psi(\lambda r)\le \psi(r)$ for  $\lambda \in (0,1)$. By defining $C_U:=a_6$ and $\beta:=2(\delta_2\vee \delta_4)$, it follows from the right-hand side inequality above that
\begin{equation}\label{eq:beta}
\psi(\lambda r)\le C_U(1+\lambda^{\beta})\psi(r)\, ,\quad \lambda >0, r>0\, .
\end{equation}
Moreover, again by \cite[(2.11)]{KSV8} we conclude that for $\lambda \le 1$ and all $r>0$
$$
\frac{f(\lambda r)}{f(r)}\ge a_6^{-1}\lambda^{\delta_2\vee \delta_4}\, ,
$$
which gives
\begin{equation}\label{eq:beta-left}
\psi(\lambda r)\ge C_U^{-1}\lambda^{\beta}\psi(r)\, ,\quad \lambda \le 1\, ,r>0\, .
\end{equation}

From now on we assume that $d\le 2$, and that $f$ is a complete Bernstein function satisfying {\bf (H1)} and {\bf (H2)}. Inequality \eqref{eq:beta} replaces \eqref{eq:square} which for $d\le 2$ is not good enough. Inequality \eqref{eq:beta-left} replaces \eqref{eq:square-left}.
It follows from \cite[Theorem 3.4(b)]{KSV8} that $X$ is transient if $d>\beta$. Moreover, the Green function exists and satisfies \eqref{lem:green_est}.
The inequality \eqref{eq:green_scaling} reads
$$
\frac{1}{2 C_U C_G^2}\lambda^{\beta-d}G(x)\leq G(\lambda x)\leq \frac{C_G^2}{C_L}\lambda^{\alpha-d}G(x)\, ,\qquad \lambda\in (0,1]\, ,\  x\in \R^d\setminus\{0\}\, .
$$
We comment now the other preliminary results from Section \ref{s:preliminaries}. Lemma \ref{l:long-jump} is valid in all dimensions (see \cite[Corollary 2]{G}), as well as Lemmas \ref{lem:maxpr} and \ref{l:equilibrium-potential}. Discussion of the scaling goes through unchanged (see \cite[Section 2]{KSV8}). The statements and proofs of Lemmas \ref{lem:scaling}--\ref{l:capacity-scaling} are the same.

We also note that the global scale invariant Harnack inequality is valid (see \cite[Theorem 3.7]{KSV8}).
Moreover, it is shown in \cite{KSV8} that the process $X$ satisfies conditions (1.4), (1.13) and (1.14) from \cite{CK}. Hence by \cite[Proposition 4.14]{CK}, nonnegative harmonic functions are H\"older continuous.

With these preparations we see that Theorems \ref{t:theorem-1}, \ref{t:theorem-2}, \ref{t:aikawa-criterion} and \ref{t:theorem-poisson} are valid as stated with minor changes in proofs --  \eqref{eq:square} is replaced with \eqref{eq:beta} giving $\beta$ instead of 2.

By considering the examples from the introduction, we see that some of them do not satisfy the current assumptions. Here are the ones where the process $X$ is not transient: (i) Brownian motion, $\psi(x)=|x|^2$; (ii) Isotropic stable process with $\psi(x)=|x|^{\beta}$ for $d\le \beta$; (iii) Sum of two independent isotropic stable processes, $\psi(x)=|x|^{\alpha}+|x|^{\beta}$, $d\le \alpha, \beta \le 2$; (iv) Truncated $\beta$-stable process. An example of a transient process not satisfying our assumptions is an independent sum of a Brownian motion and isotropic $\beta$-stable process with $\beta <d$, $\psi(x)=|x|^2+|x|^{\beta}$. In particular, the right-hand side inequality in {\bf (H1)} is not true, and consequently, \eqref{eq:beta} fails.

\section{Appendix}\label{s:appendix}
\subsection{Wiener's criterion for minimal thinness at infinity}
Let $X=(X_t,\P_x)$ be a subordinate Brownian motion in $\R^d$, $d\ge 3$, satisfying the weak lower scaling condition \eqref{eq:lower_scaling_condition}. The goal of this appendix is to outline a proof of Theorem \ref{t:aikawa-criterion}.

The concept of minimal thinness is defined for points on the minimal Martin boundary. For the exposition in the classical case of Brownian motion  we refer to \cite{AG}, while for a class of Markov processes the theory was developed in \cite{Fol}. Our first goal is
to show that the Martin boundary of $\R^d$ with respect to $X$ consists of one point that we call infinity and denote by $\infty$ (the fact tacitly used in introduction). This will follow from the next Liouville-type result. In the proof we use a global scale invariant Harnack inequality proved in \cite[Theorem 1]{G}: There exists a constant $C>0$ such that for every $r>0$ and every function $u:\R^d\to [0,\infty)$ which is harmonic in $B(0,2r)$,
$$
\sup_{x\in B(0,r)} u(x) \le C \inf_{x\in B(0,r)} u(x)\, .
$$

\begin{Lem}\label{p:liouville}
If $u:\R^d\to [0,\infty)$ is harmonic in $\R^d$ with respect to $X$, then $u$ is a constant function.
\end{Lem}
\proof Without loss of generality we may assume that $\inf_{x\in \R^d}u(x)=0$ (otherwise subtract the infimum). Then $u(x)=0$ for every $x\in \R^d$. Suppose not, and without loss of generality assume that $u(0)>0$ (note that $X$ is translation invariant). Then for every $r>0$,
$$
0<u(0)\le\sup_{x\in B(0,r)}u(x)\le C\inf_{x\in B(0,r)}u(x)\, .
$$
By letting $r\to \infty$ we obtain $0<u(0)\le 0$ -- contradiction. \qed

This shows that the Martin boundary, and hence the minimal Martin boundary, of $\R^d$ consists of a single point that we denote by $\infty$.

For $A\subset \R^d$, let $T_A:=\inf\{t> 0:\, X_t\in A\}$. A Borel set $A\subset \R^d$ is said to be \emph{minimally thin} with respect to $X$ at $\infty$ if there exists $x\in \R^d $ such that
$$
P_A 1(x):=\P_x(T_A<\infty)<1\, .
$$
Here $P_A$ denotes the hitting operator to $A$ for $X$: $P_A u(x)=\E_x[u(X_{T_A}), T_A<\infty]$.
In potential-theoretic language, $P_A 1=
\wh{R}^A_1$ - the balayage of $1$ onto $A$. Recall that if $u$ is an excessive function with respect to $X$ and $A\subset \R^d$, then the reduit of $u$ onto $A$ is defined by
$$
R^A_u=\inf\{s:\, s
  \ \text{excessive and } s\ge u \textrm{ on }A\}.
$$
Probabilistically, $R^A_u (x)=\E_x[u(X_{D_A}),D_A<\infty]$, where $D_A=\inf\{t\ge 0:\, X_t\in A\}$ is the debut of $A$. The balayage $\wh{R}^A_u$ is the lower semi-continuous regularization of the reduit $R^A_u$. It holds that $\wh{R}^A_u\le R^A_u\le u$, and $R^A_u=u$ on $A$ (see \cite{BH} for details).

\begin{Prop}\label{p:thinness}  The following are equivalent:

\noindent (a)
$A$ is minimally thin at $\infty$;

\noindent (b) There exists a potential $u=G\mu$  such that
\begin{equation}\label{e:thinness}
\liminf_{x\to \infty,\ x\in A} u(x)>0\, .
\end{equation}

\noindent (c) There exists a potential $u=G\mu$ such that
\begin{equation}\label{e:thinness2}
\liminf_{x\to \infty,\ x\in A} u(x)=+\infty \, .
\end{equation}
\end{Prop}
\proof We sketch the proof following the proof of \cite[Theorem 9.2.6]{AG}. Clearly, (c) implies (b). Assume that (b) holds. Then $\liminf_{x\to \infty, x\in A} u(x)=:a >0$. Hence, there exists a Martin topology neighborhood $W$ of $\infty$ such that $u\ge a/2$ on $A\cap W$. If $\wh{R}^{A\cap W}_1 =1$, then $u\ge \wh{R}^{A\cap W}_u\ge a/2>0$ everywhere, which implies that $a/2$ is a harmonic minorant of $u$. This is impossible, since $u$ is a potential. Hence $\wh{R}^{A\cap W}_1 \neq 1$, i.e., $A$ is minimally thin at $\infty$.

Suppose that (a) holds. By \cite[Lemma 2.7]{Fol}, there exists an open subset $U\subset \R^d$ such that $A\subset U$, and $U$ is minimally thin at $\infty$. By the analog of \cite[Theorem 9.2.5]{AG}, there is a sequence $(W_n)_{n\ge 1}$ of Martin topology open
neighborhoods of $\infty$ shrinking to $\infty$ such that $\wh{R}^{U\cap W_n}_1(0)\le 2^{-n}$. Let $u_1:=\sum_{n=1}^{\infty} \wh{R}^{U\cap W_n}_1$. Then
$u_1$ is a sum of potentials, hence a potential itself since $u_1(0)<\infty$. Further, $\wh{R}^{U\cap W_n}_1=1$ on the open set $U\cap W_n$. Therefore, $u_1(x)\to \infty$ as $x\to \infty$, $x\in U$. Thus (c) holds. \qed

Before we proceed, let us record an estimate on the Green function. Let $\eta \in (0,1)$. By Lemma \ref{lem:green_est},
$$
G(\eta|y|)\le C_G^2 \eta^{-d} \frac{\psi^{\star}(|y|^{-1})}{\psi^{\star}(\eta^{-1}|y|^{-1})}\, G(|y|)\, .
$$
By \eqref{eq:square},
$$
\psi^{\star}(|y|^{-1})=\psi^{\star}(\eta \eta^{-1}|y|^{-1})\le 2(1+\eta^2)\psi_0^{\star}(\eta^{-1}|y|^{-1})\le 4\psi^{\star}(\eta^{-1}|y|^{-1})\, .
$$
Thus, for $\eta\in (0,1)$,
\begin{equation}\label{e:green-estimate-reverse}
G(\eta|y|)\le 4C_G^2 \eta^{-d} G(|y|)\, ,\quad y\in \R^d, y\neq 0\, .
\end{equation}

The following proposition is an analog of \cite[Proposition V. 4.15]{BH}. The proof is similar and uses Lemma \ref{lem:green_est} and \eqref{e:green-estimate-reverse}. We omit the proof.
\begin{Prop}\label{p:minthin-criterion-1}
Let $E\subset \R^d$ such that $0\notin E$. Let $(s_n)_{n\ge 1}$ be a sequence in $(0,\infty)$ and $\delta \in (0,1)$ such that $s_{n+1}\ge \delta^{-1} s_n$ for every $n\ge 1$. Define
$$
A_n=E\cap \{x\in \R^d:\, s_n\le |x|<s_{n+1}\}\, ,\quad n\ge 1\, .
$$
Then $E$ is minimally thin at $\infty$ if and only if $\sum_{n=1}^{\infty}R^{A_n}_1(0)<\infty$.
\end{Prop}
\begin{Cor}\label{c:wiener-criterion}(Wiener's criterion)
Let $E\subset \R^d$ such that $0\notin E$, and let $\lambda >1$. For $n\in \N$ define $A_n=\{x\in \R^d:\, \lambda^{n+1}\le |x|<\lambda^n\}$. Then $E$ is minimally thin at $\infty$ if and only if
\begin{equation}\label{e:wiener-criterion}
\sum_{n=1}^{\infty} G(\lambda^n)\mathrm{Cap}(E\cap A_n)<\infty\, .
\end{equation}
\end{Cor}
Again, this is an analog of \cite[Corollary V. 4.17]{BH} with the same proof.

\subsection{Quasiadditivity of capacity}
In this subsection we indicate the main steps in the proof that the capacity $\mathrm{Cap}$ is quasiadditive with respect to a certain Whitney decomposition. More precisely, if $\{Q_m\}_{m\ge 1}$ is a Whitney decomposition of $\R^d\setminus\{0\}$, then there exists a constant $C_5>0$ such that for every $E\subset \R^d\setminus\{0\}$,
\begin{equation}\label{e:quasiadditivity}
\mathrm{Cap}(E)\ge C_5\sum_{m=1}^{\infty}\mathrm{Cap}(E\cap Q_m)\, .
\end{equation}

This will follow from \cite[Theorem 7.1.3]{AE} once we check conditions of that theorem. The first condition is that $\{Q_m,\wt{Q}_m\}$, where $\wt{Q}_m$ is the double cube, is a quasidisjoint decomposition of $\R^d\setminus\{0\}$ (see \cite[pages 146, 147]{AE} for details). This is clear. The second condition is that the kernel, in our case $G(x,y)$, satisfies the Harnack property with respect to $\{Q_m,\wt{Q}_m\}$, namely that $G(x,y)\asymp G(x',y)$ for $x,x'\in Q_m$, $y\in \wt{Q}_m^c$. This is also clear since distances $|x-y|$ and $|x'-y|$ are comparable.  What remains to show is that there exists a measure $\sigma$ comparable to $\mathrm{Cap}$ with respect to $\{Q_m\}$. This means that there exists a constant $C_6>0$ such that
\begin{eqnarray}
\sigma(Q_m)&\asymp &\mathrm{Cap}(Q_m)\,  \quad \textrm{for all }Q_m\, ,\label{e:def-of-comparable-measure1}\\
\sigma(E)&\le &C_6 \mathrm{Cap}(E)\, \quad \textrm{for every Borel set }E\, .\label{e:def-of-comparable-measure2}
\end{eqnarray}

Let us define the measure $m_{\psi}$ in the following way:
\begin{equation}\label{e:measure-m}
m_{\psi}(E)=\int_E \psi^{\star}(|x|^{-1})\, dx\, .
\end{equation}
Note that for any $r>0$,
$$
m_{\psi}(B(0,r))=\int_{B(0,r)}\psi^{\star}(|x|^{-1})\, dx \le c_1\int_0^r t^{d-1}\psi^{\star}(t^{-1})\, dt \le c_2 r^{d}\psi^{\star}(r^{-1})
\le c_3\frac{1}{G(r)}\, .
$$

The proof of the following statements is similar to the proof of \cite[Lemma 1]{Aik}.
\begin{Lem}\label{l:estimate-of-m}
\begin{itemize}
    \item[(i)] If $0<r<|x|/2$, then $m_{\psi}(B(x,r))\asymp r^d \psi(|x|^{-1})$.
    \item[(ii)] If $r\ge |x|/2$, then $m_{\psi}(B(x,r))\asymp \frac{1}{G(r)}\asymp r^d \psi(r^{-1})$.
    \item[(iii)] In both cases, there exists a constant $C>0$ such that $m_{\psi}(B(x,r))\le C\frac{1}{G(r)}$ for all $r>0$ and all $x\in \R^d$.
    \item[(iv)] There exists $C=C(d)>0$ such that for all $r>0$ and all $x\in \R^d$ $m_{\psi}(B(x,2r))\le C m_{\psi}(B(x,r))$, i.e.~$m_{\psi}$ is a doubling measure.
    \item[(v)] Let $p\in (1,\infty)$ and $1/p+1/q=1$. There exist a constant $C>0$ such that $\int_{B(x,r)}|y|^{-d/q}\psi(|y|^{-1})^{1/p}\, dy \le C \left(m_{\psi}(B(x,r))\right)^{1/q}$.
\end{itemize}
\end{Lem}

\begin{Rem}\label{r:comparable-measure}
Let $Q$ be a Whitney cube for $\R^d\setminus \{0\}$. This means that $0\notin Q$, and the sidelength $\ell(Q)$ of $Q$ is comparable to its diameter $\mathrm{diam}(Q)$. If $x_0$ denotes the center of $Q$, then $B(x_0, c_1 \mathrm{diam}(Q))\subset Q \subset B(x_0, c_2 \mathrm{diam}(Q))$ (for universal $0<c_1<c_2$). Hence
$$
\mathrm{Cap}(Q)\asymp \mathrm{Cap}(B(x_0,\mathrm{diam}(Q))\asymp \frac{1}{G(\mathrm{diam}(Q))}\, .
$$
On the other hand, since $|x_0|$ is comparable with $\mathrm{diam}(Q)$, by Lemma \ref{l:estimate-of-m}
$$
m_{\psi}(Q)\asymp m_{\psi}(B(x_0, \mathrm{diam}(Q))\asymp \frac{1}{G(\mathrm{diam}(Q))}\, .
$$
This shows that $m_{\psi}$ satisfies \eqref{e:def-of-comparable-measure1} with respect to the family of Whitney cubes for $\R^d\setminus \{0\}$.
\end{Rem}

Now we need an estimate for the difference of two Green functions. This is the only place where we use the assumption that $X$ is a subordinate Brownian motion.

\begin{Lem}\label{l:gradient-estimate}
There exists a constant $c>0$ such that for all $x,y\in \R^d$ we have
\begin{equation}\label{e:gradient-estimate}
\left|G(x)-G(y)\right| \le c \left(\frac12 \wedge \frac{|x-y|}{|x|\wedge|y|}\right)G\left(|x|\wedge|y|\right)\, .
\end{equation}
\end{Lem}
\proof Let $U(dt)$ denote the potential (renewal) measure of the subordinator. Then (see, e.g.~\cite[(5.47)]{SV})
$$
G(x)=\int_{[0,\infty)}(4\pi t)^{-d/2}e^{-\frac{|x|^2}{4t}}\, U(dt)\, , \quad x\in \R^d\, .
$$
We look separately at two cases.

Case 1: $|x-y|\ge \frac12\left(|x|\wedge |y|\right)$. Then by monotonicity of $G$ and Lemma 2.1 (2.6), we trivially have that
$$
\left|G(x)-G(y)\right|\le G(|x|)+G(|y|)\le 2 G\left(\frac12(|x|\wedge |y|)\right)\le c_1 G(|x|\wedge|y|)\, .
$$

Case 2: $|x-y|\le \frac12\left(|x|\wedge |y|\right)$. By the mean value theorem we have that
$$
e^{-\frac{|x|^2}{4t}}-e^{-\frac{|y|^2}{4t}}=- \frac{2\langle x+\vartheta(y-x), y-x\rangle}{4t}\, e^{-\frac{|x+\vartheta(y-x)|^2}{4t}}\, ,
$$
where $\vartheta=\vartheta(x,y,t)\in (0,1)$. We estimate the absolute value above as follows:
\begin{eqnarray*}
\left|e^{-\frac{|x|^2}{4t}}-e^{-\frac{|y|^2}{4t}}\right|&\le &|y-x|\, \frac{| x+\vartheta(y-x)|}{2t}\, e^{-\frac{|x+\vartheta(y-x)|^2}{4t}}\\
&\le & \frac{2|y-x|}{| x+\vartheta(y-x)|} \, \frac{| x+\vartheta(y-x)|^2}{4t}\, e^{-\frac{|x+\vartheta(y-x)|^2}{4t}}\\
&\le & \frac{2|y-x|}{| x+\vartheta(y-x)|}\, e^{-\frac{|x+\vartheta(y-x)|^2}{8t}}\, ,
\end{eqnarray*}
where in the last line we used the elementary inequality $s e^{-s}\le e^{-s/2}$, $s>0$. Since
$$
|x+\vartheta(y-x)|\ge |x|-\vartheta|y-x|\ge |x|-\frac{|x|\wedge |y|}{2}\ge \frac{|x|\wedge |y|}{2}\, ,
$$
we get that
$$
\left|e^{-\frac{|x|^2}{4t}}-e^{-\frac{|y|^2}{4t}}\right|\le \frac{4|y-x|}{|x|\wedge|y|}\, e^{-\frac{\left(\frac{|x|\wedge|y|}{2\sqrt{2}}\right)^2}{4t}}\, .
$$
Therefore,
\begin{eqnarray*}
\left|G(x)-G(y)\right|&\le &\int_{[0,\infty)}(4\pi t)^{-d/2} \left|e^{-\frac{|x|^2}{4t}}-e^{-\frac{|y|^2}{4t}}\right|\, U(dt)\\
&\le & \frac{4|y-x|}{|x|\wedge|y|}\, \int_{[0,\infty)} (4\pi t)^{-d/2}e^{-\frac{\left(\frac{|x|\wedge|y|}{2\sqrt{2}}\right)^2}{4t}}\, U(dt) \\
&=&  \frac{4|y-x|}{|x|\wedge|y|}\,  G\left(\frac{|x|\wedge|y|}{2\sqrt{2}}\right) \le c_2 \frac{|y-x|}{|x|\wedge|y|}\, G(|x|\wedge  |y|)\, ,
\end{eqnarray*}
where in the last inequality we have again used (2.6). \qed

Let us define the maximal function with respect to the measure $m_{\psi}$. For $f:\R^d\to \R$, let
\begin{equation}\label{e:maximal-function}
\mathcal{M}_{\psi}f(x)=\sup \frac{1}{m_{\psi}(Q)}\int_Q f(y)\, m_{\psi}(dy)\, ,
\end{equation}
where the supremum is taken over all cubes $Q$ containing $x$.

From now on, given a cube $Q$, we denote by $\wt{Q}$ the double cube, that is the cube centered at the same point as $Q$ but the size twice long.
\begin{Lem}\label{l:estimate-by-max-function}
\begin{itemize}
\item[(i)] There exists a constant $C>0$ such that for all $f:\R^d \to \R$, all $r>0$ and all $x_0\in \R^d$,
\begin{equation}\label{e:estimate-by-max-function-1}
\sup_{x,x'\in B(x_0,r)}\int_{\R^d\setminus B(x_0,2r)} |G(x,y)-G(x',y)| f(y)\, m_{\psi}(dy)\le C M_{\psi}f(\tilde{x}_0)\, ,
\end{equation}
where $\wt{x}_0$ is any point in $B(x_0,r)$.
\item[(ii)] There exists a constant $C>0$ such that for all $f:\R^d \to \R$ and all cubes $Q$
\begin{equation}\label{e:estimate-by-max-function-2}
\sup_{x,x'\in Q}\int_{\R^d\setminus \wt{Q}} |G(x,y)-G(x',y)| f(y)\, m_{\psi}(dy)\le C M_{\psi}f(\wt{x}_0)\, ,
\end{equation}
where $\wt{x}_0$ is any point in $B(x_0,r)$.
\item[(iii)] There exists a constant $C>0$ such that for all cubes $Q$
\begin{equation}\label{e:estimate-by-max-function-2}
\sup_{x,x'\in Q}\int_{\R^d\setminus \wt{Q}} |G(x,y)-G(x',y)| \, m_{\psi}(dy)\le C \, .
\end{equation}
\end{itemize}
\end{Lem}
The proof of (i) and (ii) follows the method of the proof of \cite[Lemma 14.4]{AE}, while (iii) follows by taking $f\equiv 1$.

Introduce the operator
\begin{equation}\label{e:operator-T}
T_{\psi}f(x):=\int_{\R^d}G(x,y)f(y)\, m_{\psi}(dy)\, ,
\end{equation}
and the norms
$$
\|f\|_{p,\psi}:=\left(\int_{\R^d}|f|^p\, dm_{\psi}\right)^{1/p}\, .
$$
The goal is to prove the following theorem.

\begin{Thm}\label{t:T-bounded}
\begin{itemize}
\item[(i)] Let $1<p<\infty$. There exists a constant $C>0$ such that $\|T_{\psi}f\|_{p,\psi}\le C \|f\|_{p,\psi}$.
\item[(ii)] There exists a constant $C>0$ such that for every $\lambda >0$, and every measure $\mu$,
\begin{equation}\label{e:weak-1-1}
\lambda  m_{\psi}\left(\left\{x\in \R^d:\, \int_{\R^d}G(x,y)\, \mu(dy) >\lambda\right\}\right) \le C \mu(\R^d)\, .
\end{equation}
\end{itemize}
\end{Thm}
Part (ii) of Theorem \ref{t:T-bounded} immediately implies that condition \eqref{e:def-of-comparable-measure2} is satisfied for the measure $m_{\psi}$. Indeed, suppose that $\mu$ is a measure such that $G\mu \ge 1$ on $E$. With $\lambda=1$, \eqref{e:weak-1-1} implies  $$
m_{\psi}(E)\le m_{\psi}(G\mu \ge 1)\le C\mu(\R^d)\, .
$$
Since $\mathrm{Cap}(E)=\inf\{\mu(\R^d):\, G\mu \ge 1 \textrm{ on } E\}$, we get that $m_{\psi}(E)\le C\mathrm{Cap}(E)$.

Theorem \ref{t:T-bounded} is proved through a series of lemmas that exactly follow Lemmas 3, 5 and 6 in \cite{Aik}. Recall that $\alpha>0$ is the exponent form the weak lower scaling condition \eqref{eq:lower_scaling_condition}. Without loss of generality we may assume that $0<\alpha<2$.

\begin{Lem}\label{l:lemma1-for-theorem}
Suppose that $p\in (1,\infty)$ satisfies $1<p<\frac{d}{d-\alpha}$, and let $\frac1p+\frac1q=1$. There exists a constant $C>0$ such that for  any cube $Q$ in $\R^d$,
\begin{equation}\label{e:lemma1-for-theorem}
\int_Q|T_{\psi}f(x)|\, m_{\psi}(dx) \le C \|f\|_{q,\psi}\, m_{\psi}(Q)^{1/p}\, .
\end{equation}
\end{Lem}

Note that if $\gamma=(1+p)^{-1}$, then $p=\frac{1}{\gamma}-1$, and $1<p<\frac{d}{d-\alpha}$ is equivalent to $\frac{d-\alpha}{2d-\alpha}<\gamma<\frac12$.

\begin{Lem}\label{l:lemma2-for-theorem}
Let $\frac{d-\alpha}{2d-\alpha}<\gamma<\frac12$. There exists a constant $C>0$ such that if $0<\epsilon <1$, $Q$ is a cube and $\wt{Q}$ the double cube, $f\ge 0$ with $\mathrm{supp}(f)\subset \wt{Q}$ and $\|f\|_{1,\psi}\le \epsilon m_{\psi}(\wt{Q})$, then
$$
m_{\psi}\left(\left\{x\in Q:\, T_{\psi}f(x)>1\right\}\right)\le C e^{1-\gamma}m_{\psi}(Q)\, .
$$
\end{Lem}

\begin{Lem}\label{l:lemma3-for-theorem}
Let $\frac{d-\alpha}{2d-\alpha}<\gamma<\frac12$. There exist constants $B>0$ and $C>0$ such that if $\lambda >0$, $0<\epsilon <1$, $f\ge 0$ and a cube $Q$ has a point $x'$ such that $T_{\psi}f(x')\le \lambda$, then
$$
m_{\psi}\left(\left\{x\in Q:\, T_{\psi} f(x)>B\lambda, \mathcal{M}_{\psi}f(x)\le \epsilon \lambda\right\}\right) \le C\epsilon^{1-\gamma} m_{\psi}(Q)\, .
$$
\end{Lem}

With Lemmas \ref{l:lemma1-for-theorem}--\ref{l:lemma3-for-theorem}, the proof of Theorem \ref{t:T-bounded} is the same as the proof of \cite[Theorem 3]{Aik}.

\subsection{Aikawa's Wiener-type criterion for minimal thinness at infinity}
Let $\{Q_m\}_{m\ge 1}$ be a Whitney decomposition of $\R^d\setminus\{0\}$ with cubes of size $3^n$, $n\in \Z$. By the previous subsection, there exists a constant $C>1$ such that for every $E\subset \R^d$,
\begin{equation}\label{e:quasi-add}
C^{-1}\sum_{m=1}^{\infty}\mathrm{Cap}(E\cap Q_m)\le \mathrm{Cap}(E)\le C\sum_{m=1}^{\infty}\mathrm{Cap}(E\cap Q_m)\, .
\end{equation}

\emph{Proof of Theorem \ref{t:aikawa-criterion}.} Take $\lambda=3$ in Corollary \ref{c:wiener-criterion} and define $A_n=\{x\in \R^d:\, 3^n\le |x|<3^{n+1}\}$. If $A_n\cap Q_m\neq \emptyset$, then $\mathrm{diam}(Q_m)\asymp 3^n$, and therefore $G(3^n)\asymp G(\mathrm{diam}(Q_m))$. Thus,
\begin{eqnarray*}
\sum_n G(3^n)\mathrm{Cap}(E\cap A_n)&\asymp &\sum_nG(3^n)\sum_m \mathrm{Cap}(E\cap A_n\cap Q_m)\\
&=&\sum_m \sum_n G(3^n)\mathrm{Cap}(E\cap A_n\cap Q_m)\\
&\asymp &\sum_m \sum_{n,A_n\cap Q_m\neq \emptyset} G(\mathrm{diam}(Q_m))\mathrm{Cap}(E\cap A_n\cap Q_m)\\
&= &\sum_m  G(\mathrm{diam}(Q_m))\sum_{n,A_n\cap Q_m\neq \emptyset}\mathrm{Cap}(E\cap A_n\cap Q_m)\\
&\asymp &\sum_m  G(\mathrm{diam}(Q_m)) \mathrm{Cap}(E\cap Q_m)\, .
\end{eqnarray*}
The first line follows from \eqref{e:quasi-add}. One inequality in the last line is subadditivity of capacity. For another we argue as follows. There exists $N\in \N$ such that for every $Q_m$, $\sum_{n, A_n\cap Q_m\neq \emptyset}1=\sum_n 1_{A_n\cap Q_m}\le N$. Hence,
$\sum_{n, A_n\cap Q_j\neq \emptyset}\mathrm{Cap}(E\cap A_n\cap Q_j) \le N \mathrm{Cap}(E \cap Q_j)$.
\qed

{\bf Acknowledgement:}
We thank Stephen Gardiner for explaining the details of  scaling in \cite[Section 5]{GG}.
We are also grateful to the referee for valuable comments.

\end{document}